% \magnification=1200

% qualche utile definizione

\input epsf.tex

\def\2{{1\over 2}}

\def\a{\alpha}
\def\b{\beta}
\def\g{\gamma}

\def\s{\sigma}
\def\e{\epsilon}
\def\l{\lambda}

\def\fun#1#2#3{#1\colon #2\rightarrow #3}

\def\frac#1#2{{{#1} \over {#2}}}

\def\st{\;\colon\;}
\def\tends{\rightarrow}

\def\dx{\hbox{{\rm d}$x$}}

\def\dr{ {\rm d} }

\def\R{{\bf R}}
\def\N{{\bf N}}
\def\Z{{\bf Z}}

\def\T{{\bf T}}

\def\W{{\cal W}}

\def\thm#1{\vskip 1 pc\noindent{\bf Theorem #1.\quad}\sl}
\def\lem#1{\vskip 1 pc\noindent{\bf Lemma #1.\quad}\sl}
\def\prop#1{\vskip 1 pc\noindent{\bf Proposition #1.\quad}\sl}

\def\proof{\rm\vskip 1 pc\noindent{\bf Proof.\quad}}
\def\fin{\par\hfill $\backslash\backslash\backslash$\vskip 1 pc}
\def\txt#1{\quad\hbox{#1}\quad}

\def\L{{\cal L}}

\def\G{\Gamma}
\def\s{\sigma}

\def\2{\frac{1}{2}}
\def\inn#1#2{{\langle #1 ,#2\rangle}}

\def\part{{\partial_{x}}}
\def\div{{\rm div}}

\def\ac{{\cal A}}
\def\fc{{\cal F}}
\def\uc{{\cal U}}
\def\vc{{\cal V}}
\def\Pc{{\cal P}(\T^d)}
\def\Pcc{{\cal P}_2(\R^d)}
\def\si{\sigma^{(t,\psi)}}

%\def\exp{{\rm exp}}
%\def\freq{{\frac{\bar k}{\bar T}}}

% pagina di titolo 

\baselineskip= 17.2pt plus 0.6pt
% \hsize=14truecm
% \hoffset=2truecm
% \vsize=8truein
% \voffset=0.5truein
\font\titlefont=cmr17
\centerline{\titlefont Existence of solutions of the master equation}
\vskip 1 pc
\centerline{\titlefont in the smooth case}
%\vskip 1 pc
% \footnote{}{\rm Supported by 
% Italian Ministry of Education}}
\vskip 4pc
\font\titlefont=cmr12
\centerline{         \titlefont {Ugo Bessi}\footnote*{{\rm 
Dipartimento di Matematica, Universit\`a\ Roma Tre, Largo S. 
Leonardo Murialdo, 00146 Roma, Italy.}}   }{}\footnote{}{
{{\tt email:} {\tt bessi@matrm3.mat.uniroma3.it} Work partially supported by the PRIN2009 grant "Critical Point Theory and Perturbative Methods for Nonlinear Differential Equations}} 
\vskip 0.5 pc
%\centerline{{\tt email:} {\tt bessi@matrm3.mat.uniroma3.it}}
 
\par
\vskip 2pc
\centerline{\bf Abstract}
We give a different proof of a theorem of W. Gangbo and A. Swiech on the short time existence of solutions of the master equation.

\vskip 2 pc
%\magnification=1200
\centerline{\bf  Introduction}
\vskip 1 pc

Mean Field Games are games with a continuum of players, each of which sees only the "mean field" generated by the other ones. They attracted the attention of a wider set of analysts after the lectures of P. L. Lions at the Coll\`ege de France, which are available in video streaming (see also the written presentation [11]). They can model a wide array of phenomena in physics and  mathematical economics; we dwell a little on one aspect of the latter. Actually, the idea of considering a continuum of players came up naturally in mathematical economy, where it was used ([6], see also [14] for a more elementary presentation) to model the formation of prices in a market with perfect concurrence. Quoting from [6], "the essential idea of this notion is that the economy under consideration has a "very large" number of participants, and that the influence of each participant is "negligible"". 

To be more precise, let us look at the situation of [15]: we have a probability measure $\mu_s$ on the $d$-dimensional torus 
$\T^d=\frac{\R^d}{\Z^d}$ which models the distribution of the players at time $s$; we fix an initial time $t<0$, an initial distribution $\bar\mu$ and we suppose that $\mu_s$ evolves according to the continuity equation, forward in time, 
$$\left\{
\eqalign{\partial_s\mu_s+\div(X\mu_s)&=0
\quad s>t\cr
\mu_t&=\bar\mu  
}
\right.     \eqno (1)$$
where the vector field $X$ is a control which we are free to choose in the following. 

Let us call ${\cal P}(\T^d)$ the space of the Borel probability measures on $\T^d$, and let us suppose that we are given two potentials $\fun{{\cal F},{\cal U}_0}{{\cal P}(\T^d)}{\R}$. We would like the whole society to minimize the value function
$$\vc(t,\bar\mu)\colon=
\inf
\left\{
\int_t^0\dr s\left[
\int_{\T^d}\2|X^2(s,x)|^2\dr\mu_s(x)-\fc(\mu_s)
\right]   
+{\cal U}_0(\mu_0)
\right\}\eqno (2)$$
where the $\inf$ is over all curves which satisfy (1) and all controls $X$. It turns out that under suitable hypotheses on $\fc$ and 
$\uc_0$ the $\inf$ is a minimum: there is a vector field $X$ minimizing in (2); by (1), we also have a minimal trajectory $\mu_s$. 

In (2), we minimize the cost for the whole society, but what about its members? One possible notion is that of Nash equilibrium: roughly, we are on a Nash equilibrium if no one can get a better deal by a unilateral change of strategy. It happens that, in our case, the optimum for the whole society is a Nash equilibrium. Actually, under suitable hypotheses on ${\cal F}$ and 
${\cal U}_0$, we shall be able to define two functions $F(x,\mu)$ and $u_0(x,\mu)$ which, heuristically, are the "mean field" potentials felt by the particle placed at $x$, provided the other ones are distributed as $\mu$. We shall see that the drift $X$ in (1) optimal for the whole group is also best for the single particle; namely, $X(s,q)=-\partial_xv(s,q)$ where $v$ solves the Hamilton-Jacobi equation with time reversed
$$\left\{
\eqalign{
-\partial_tv(s,q)+\2|\partial_q v(s,q)|^2+F(q,\mu_s)&=0
\quad s\le 0\cr
v(0,q)&=u_0(q,\mu_0)  .
}      
\right.     \eqno (3)$$
Equivalently, the particle initially placed at $q$ minimizes its cost:
$$\int_t^0\2[|\dot q(s)|^2+F(q(s),\mu_s)]\dr s+u_0(q(0),\mu_0)    
$$
if it follows the vector field $X$. 

Since the value function $\vc(t,\mu)$ of (2) is defined on the metric space ${\cal P}(\T^d)$, this approach calls for a study of the Hamilton-Jacobi equation in metric spaces; we refer the reader to [3], [16] and [20] for three definitions of viscosity solutions of H-J in metric spaces. 

In this framework, the task is to solve the coupled equations (1) and (3); it turns out that, formally, these two equations are equivalent to the so-called master equation, i. e. formula (6) below. Heuristically, the solution of the master equation is a value function both for the single particle and the whole community. In [15] it is shown that, under suitable hypotheses on $\fc$ and 
$\uc$, the master equation has a smooth solution for $t$ negative and small and that the master equation is equivalent (this time rigorously) to (1) and (3).  

In this paper, we want to give a different proof of the results of [15]. Instead of working in $\Pc$, we take up a suggestion of [11] (see also [18], [19]) and work in the space of $L^2$ parametrizations of particles: a parametrization for $\mu$ will be a function 
$\s\in L^2([0,1)^d,\R^d)$ whose law, when projected on $\T^d$, is $\mu$. In other words, we are choosing $[0,1)^d$ as parameter space. 

We shall see that this approach is equivalent to that of [15]; as in [15], the implicit function theorem is at the core of our proof, but we are going to use it in a way that is closer to the original approach of [10]. 

We set $M=L^2([0,1)^d,\R^d)$ and denote by $AC([a,b],X)$ the set of the absolutely continuous functions from $[a,b]$ to a space $X$; throughout the paper, we shall denote by $\nabla$, $D$ and 
$d$ the gradients of functions on $\T^d$, $M$ and $\Pc$ respectively. 

We want to prove the following.

\thm{1} Let $\fun{\hat\fc,\hat\uc_0}{M}{\R}$ be respectively a potential and a final condition satisfying the hypotheses of section 2 below. Then, the following points hold. 

\noindent $1$) There is $T>0$ such that, if $t\in[-T,0]$ and 
$\psi\in M$, the minimum
$$\hat\uc(t,\psi)\colon=
\min\left\{
\int_t^0[
\2||\dot\s_s||_M^2-\hat\fc(\s_s)
]\dr s+\hat\uc_0(\s_0)\st
\s\in AC([t,0],M),\quad \s_t=\psi
\right\} \eqno (4)$$
is attained on a unique curve $\s^{(t,\psi)}\in AC([t,0],M)$. 

\noindent $2$) The maps $\fun{}{(t,\psi)}{\si}$ and 
$\fun{}{(t,\psi)}{\hat\uc(t,\psi)}$ are of class $C^2$; moreover, they are $L^2_\Z$ and $H$-equivariant in the last variable for the groups $L^2_\Z$ and $H$ defined in section 1 below. 

\noindent $3$) There are two functions of class $C^3$
$$\fun{\hat F,\hat u_0}{\T^d\times M}{\R}$$
such that, if we set
$$u(t,x,\psi)=\min\Big\{
\int_t^0[
\2|\dot q(s)|^2-\hat F(q(s),\s_s^{(t,\psi)})
]  \dr s+\hat u_0(q(0),\s_0^{(t,\psi)})\st$$
$$q\in AC([t,0],\T^d),\quad q(t)=x
\Big\}   \eqno(5)$$
then $u$ is of class $C^2$ in $[-T,0]\times\T^d\times M$ and satisfies the master equation
$$-\partial_tu(t,q,\psi)+\2|\nabla u(t,q,\psi)|^2+F(q,\psi)+
\inn{\nabla u(t,\psi(\cdot),\psi)}{D u(t,q,\psi)}_M=0   
\qquad\forall (t,x,\psi)\in[-T,0]\times\T^d\times M      \eqno (6)$$
where $\inn{\cdot}{\cdot}_M$ denotes the inner product in $M$. To districate the inner product above, we note that 

\noindent $D u(t,q,\psi)\in M$ because it is the gradient with respect to the $M$ variable; moreover, 
$\fun{}{x}{\nabla u(t,\psi(x),\psi)}$ belongs to $M$ since it is the $C^2$ function $u(t,\cdot,\psi)$ composed with $\psi$. The function $u$ is $\Z^d$-equivariant in the second variable and $L^2_\Z$ and $H$-equivariant in the last one. 

\noindent $4$) Let the law of $\psi$ be absolutely continuous with respect to the Lebesgue measure; then, for $s\in[-T,0]$ the law of 
$\s_s^{(t,\psi)}$ is absolutely continuous too. 

\noindent $5$) For $\L^d$ a. e. $x\in[0,1)^d$ we have that, for all 
$s\in[-T,0]$,
$$\dot\s_s^{(t,x)}(x)=-\nabla u(s,\s^{(t,x)}_s(x),\s^{(t,x)}_s)
  .  $$
In other words, the orbit $q(s)$ minimal in (5) coincides with 
$\s_s^{(t,\psi)}(x)$ if they start at the same point of $\T^d$; equivalently, $\fun{}{s}{\s^{(t,\psi)}_s}(x)$ minimizes the one-particle problem (5) for $\L^d$ a. e. $x\in[0,1)^d$. 

\rm

\vskip 1pc

Recently the master equation has been studied extensively, expecially from the stochastic viewpoint; we refer the reader to [7], [8], [9], [12] and [13]. 

The paper is organized as follows: section 1 contains the notation and a theorem of [11] about the relationship between differentiability on parametrizations and on measures; section 2 recalls the hypotheses used in [15] from section 6 onwards; in section 3 we recall the method of [10] for the minimum of (4), in section 4 we deal with the master equation (6).

\vskip 2pc
\centerline{\bf \S 1}
\centerline{\bf Preliminaries and notation}
\vskip 1pc

We denote by $\fun{\pi}{\R^d}{\T^d\colon=\frac{\R^d}{\Z^d}}$ the natural projection, and by 
$|\cdot |_{\T^d}$ the distance on $\T^d$ given by 
$$|x-y|_{\T^d}=
\min\{
|\tilde x-\tilde y|\st\pi(\tilde x)=x,\quad \pi(\tilde y)=y
\}  .  $$
We let $\Pc$ be the space of Borel probability measures on 
$\T^d$; if $\mu_1,\mu_2\in\Pc$, we denote by $\G(\mu_1,\mu_2)$ the set of all the Borel probability measures on $\T^d\times\T^d$ whose first and second marginals are, respectively, $\mu_1$ and 
$\mu_2$. For $\l\ge 1$ we define the $\l$-Wasserstein distance on 
$\Pc$ by 
$$\W_\l(\mu_1,\mu_2)^\l=
\min_{\g\in\G(\mu_1,\mu_2)}
\int_{\T^d\times\T^d}|x-y|_{\T^d}^\l\dr\g(x,y)  .  \eqno (1.1)$$
We refer the reader to [4] or [23] for the proof that the minimum  is attained and that $(\Pc,\W_\l)$ is a compact metric space. 

When $\l=2$ (which is the only case we consider in this paper) we denote by $\G_o(\mu_1,\mu_2)$ the set of the minimizers in (1.1). 

We want to parametrize $\mu\in\Pc$ with a map 
$\s\in M\colon= L^2([0,1)^d,\R^d)$. To do this, we begin to define 
$\Pcc$ as the set of the Borel probability measures on $\R^d$ with finite second moment. Following [19], we push forward 
$\mu\in\Pcc$ to $\tilde\mu\colon=\pi_\sharp \mu\in\Pc$. By the definition of push-forward, this is tantamount to 
$$\int_{\T^d}f(x)\dr\tilde\mu(x)=
\int_{\R^d}f(x)\dr\mu(x)\qquad\forall f\in C(\T^d,\R)  $$
where we have identified $f$ with its lift to a periodic function on 
$\R^d$. 

If $\pi_\sharp\mu_1=\pi_\sharp\mu_2=\tilde\mu$, we say with [19] that $\mu_1$ and $\mu_2$ are two representatives of 
$\tilde\mu$. By lemma 1.2 of [19], it is possible to lift any couple of measures on $\T^d$ to measures on $\R^d$ in such a way to preserve the 2-Wasserstein distance. More precisely,  if 
$\tilde\mu_1,\tilde\mu_2\in\Pc$, then there are two representatives $\mu_1,\mu_2\in\Pcc$ such that 
$\mu_1$ is supported in $[0,1]^d$, $\mu_2$ in $[-1,2]^d$ and
$$\W_2(\tilde\mu_1,\tilde\mu_2)^2=
W_2(\mu_1,\mu_2)^2\colon=
\min_{\g\in\G(\mu_1,\mu_2)}
\int_{\R^d\times\R^d}|x-y|^2\dr\g(x,y)    \eqno (1.2)$$
where we have denoted by $W_2$ the 2-Wasserstein distance on 
$\Pcc$.   

Let $\L^d$ denote the $d$-dimensional Lebesgue measure on 
$[0,1)^d$ and let $\mu\in\Pcc$; it is standard ([4] or [23]) that there is a map 
$\psi\in M$ (actually, $\psi$ is the gradient of a convex function) such that $\psi_\sharp\L^d=\mu$. The trivial converse is that, if 
$\psi\in M$, then $\psi_\sharp\L^d\in\Pcc$. The map $\psi$ is called the Brenier map, or the parametrization of $\mu$. 

For completeness' sake, we give a well-known extension of lemma 6.4 of [11].

\lem{1.1} 1) Let $\mu_1,\mu_2\in\Pcc$, let $\psi_1,\psi_2\in M$ be two parametrizations of $\mu_1$, $\mu_2$ respectively and let 
$\g\in\G(\mu_1,\mu_2)$.  Then, there is a sequence of  invertible, measure-preserving maps 
$\fun{h_n}{[0,1)^d}{[0,1)^d}$ such that 
$(\psi_1\circ h_n,\psi_2)_\sharp\L^d$ converges weak$\ast$ to 
$\g$. Moreover, for all functions $f\in C(\T^d\times\R^d,\R)$ such that $\frac{f(x,v)}{1+|v|^2}$ is bounded, we have that 
$$\int_{\R^d\times\R^d}f(x,x-y)\dr\g(x,y)=
\lim_{n\tends+\infty}\int_{[0,1)^d}
f( \psi_1\circ h_n(x),\psi_2(x)-\psi_1\circ h_n(x) )\dr x  .  \eqno (1.3)$$

\noindent 2) Let $\tilde\mu_1,\tilde\mu_2\in\Pc$ and let 
$\mu_1,\mu_2\in\Pcc$ be two representatives such that (1.2) holds. Let $\psi_1,\psi_2\in M$ be as in point 1). Then,  
$$\W_2(\tilde\mu_1,\tilde\mu_2)^2=
W_2(\mu_1,\mu_2)^2=
\inf\int_{[0,1)^d}|\psi_1\circ h(x)-\psi_2(x)|^2\dr x
\eqno (1.4)$$
where the $\inf$ is over all invertible, measure-preserving maps 
$\fun{h}{[0,1)^d}{[0,1)^d}$. 

\proof As for (1.4), the first equality comes from (1.2). For the second one, we note that, since 
$(\psi_1\circ h,\psi_2)_\sharp\L^d\in\G(\mu_1,\mu_2)$, we have that 
$$W_2(\mu_1,\mu_2)^2\le
\inf_h\int_{[0,1)^d}|\psi_2(x)-\psi_1\circ h(x)|^2\dr x  .  $$
The opposite inequality follows immediately from point 1), which we prove it in the steps below using a variation of the technique of [11].   

\noindent {\bf Step 1.} We begin to suppose that $\mu_1$ and 
$\mu_2$ are supported in a common cube, say 
$\tilde Q^l=[-l,l)^d$.  We partition $\tilde Q^l$ into smaller cubes 
$$Q_k=\frac{2kl}{2^n}+\frac{1}{2^n}\tilde Q^l$$
with $k=(k_1,\dots,k_d)\in\Z^d$ such that 
$-2^n+1\le k_i\le  2^n-1$. Next, we relabel the $Q_k$ to 
$Q_i$, with $i$ in a finite set of $\N$. 

In the step 3, 4 and 5 below we are going to find maps $h_n$ such that
$$\L^d[(\psi_1\circ h_n,\psi_2)^{-1}(Q_i\times Q_j)]=
\g(Q_i\times Q_j)\txt{for all} i,j  .  \eqno (1.5)$$
Using the fact that the sides of $Q_i$ have length $\frac{2l}{2^n}$ and that $\mu_1$ and $\mu_2$ are supported in $\tilde Q_l$, the formula above easily implies that
$(\psi_1\circ h_n,\psi_2)_\sharp\L^d$ converges to $\g$ in the weak$\ast$ topology. Formula (1.3) now follows because $\g$ and 
$(\psi_1\circ h_n,\psi_2)_\sharp\L^d$ are supported in 
$\tilde Q^l\times\tilde Q^l$, a compact set on which 
$\fun{}{(x,y)}{f(x,y-x)}$ is continuous. 

\noindent{\bf Step 2.} Before showing (1.5) for the case with bounded support, let us show how it implies (1.3) in the general case. 

Let $\fun{h}{[0,1)^d}{[0,1)^d}$ be measure preserving. The equality below comes from the definition of push-forward; in the inequality, $\tilde Q^l$ is the cube of step 1.
$$\left\vert
\int_{[0,1)^d} f(\psi_1\circ h(x),\psi_2(x)-\psi_1\circ h(x))\dr x-
\int_{\R^d\times\R^d} f(x,y-x)\dr\g(x,y)
\right\vert   =  $$
$$\left\vert
\int_{\R^d\times\R^d} f(x,y-x)\dr(\psi_1\circ h,\psi_2)_\sharp\L^d(x,y)-
\int_{\R^d\times\R^d} f(x,y-x)\dr\g(x,y)
\right\vert   \le   $$
$$\int_{(\tilde Q^l\times\tilde Q^l)^c}
|f(x,y-x)|\dr(\psi_1\circ h,\psi_2)_\sharp\L^d(x,y)+  
\eqno (1.6)_a$$
$$\int_{(\tilde Q^l\times\tilde Q^l)^c}
|f(x,y-x)|\dr\g(x,y)+   \eqno (1.6)_b$$
$$\left\vert
\int_{(\tilde Q^l\times\tilde Q^l)}
f(x,y)\dr(\psi_1\circ h,\psi_2)_\sharp\L^p(x,y)-
\int_{(\tilde Q^l\times\tilde Q^l)}
f(x,y-x)\dr\g(x,y)
\right\vert   .  \eqno (1.6)_c$$
Let $\e>0$; from the formula above we see that (1.3) follows if we prove that we can find $l\in\N$ such that
$$(1.6)_a<\e  $$
for all measure-preserving $h$, 
$$(1.6)_b\le\e$$
and that, once $l$ is fixed in this way, we can find a measure-preserving $h$ such that
$$(1.6)_c\le\e  .  $$
The last formula comes immediately from step 1; $(1.6)_b<\e$ follows because the measure $|f(x,y-x)|\g$ is finite and 
$\cap_l(\tilde Q_l\times\tilde Q_l)^c=\emptyset$. 

As for $(1.6)_a\le\e$, it suffices to prove that 
$|f(x,y-x)|(\psi_1\circ h,\psi_2)_\sharp\L^d$ is a tight set of measures as $h$ varies in the measure-preserving maps of 
$[0,1)^d$. By our hypotheses on $f$, this follows if we show that 
$(1+|y-x|^2)(\psi_1\circ h,\psi_2)_\sharp\L^d$ is tight. This is equivalent to say that $|\psi_1\circ h-\psi_2|^2$ is uniformly integrable as $h$ varies among the measure-preserving maps, which follows if we prove that $|\psi_1\circ h|^2$ is uniformly integrable; we leave the easy proof of this to the reader. 

\noindent{\bf Step 3.} In this step, we define the pre-images of the cubes $Q_i$, which the map $h_n$ of step 1 will permute in a Rubik cube fashion. We set 
$$A_i=\psi_1^{-1}(Q_i)\subset[0,1)^d,\qquad
B_i=\psi_2^{-1}(Q_i)\subset[0,1)^d  .  $$
The equalities on the left in the two formulas below follow since 
$\g\in\G(\mu_1,\mu_2)$; those on the right come from the fact that $\mu_j=(\psi_j)_\sharp\L^d$ for $j=1,2$.
$$\g(Q_i\times[-l,l)^d)=\mu_1(Q_i)=\L^d(A_i),\qquad
\g([-l,l)^d\times Q_i)=\mu_2(Q_i)=\L^d(B_i)  .  \eqno (1.7)$$

In the next two steps, we shall settle the first row of cubes, say 
$\{ A_i\times B_1 \}_i$. The idea is to partition $B_1$ into sets 
$B_{i,1}$ and to find sets $A_{i,1}\subset A_i$ such that 
$\L^d(A_{i,1})=\L^d(B_{i,1})=\g(Q_i\times Q_1)$; then, we shall send $A_{i,1}$ into $B_{i,1}$ by a measure-preserving map. We shall see that this yields (1.5) for $j=1$. 

\noindent{\bf Step 4.} We assert that we can find  sets 
$A_{i,1}\subset A_i$  such that
$$\L^d(A_{i,1})=\g(Q_i\times Q_1)
\txt{and} 
\sum_i\L^d(A_{i,1})=\L^d(B_1)  .  \eqno (1.8)$$
Note that the sets $A_{i,1}$ are disjoint since the $A_i$ are disjoint. Moreover, we can find sets $B_{i,1}\subset B_1$ such that 
$$\left\{
\matrix{
\L^d(B_{i,1})=\L^d(A_{i,1})\cr
\txt{the $B_{i,1}$ are disjoint}\cr
\L^d\left( B_1\setminus\bigcup_i B_{i,1}\right) =0\cr
B_{i,1}\supset A_{i,1}\cap B_1\cr
B_{i,1}\cap A_{j,1}=\emptyset\txt{if}j\not=i  .
}
\right.   \eqno (1.9)$$
We begin to show that the first equality of (1.8) implies the second one: the first equality below follows since the $Q_i$ partition 
$[-l,l)^d$, the second one follows since $\g$ has $\mu_2$ as the second marginal, the third one since $(\psi_2)_\sharp\L^d=\mu_2$ and the fourth one from the definition of $B_1$. 
$$\sum_i\g(Q_i\times Q_1)=
\g([-l,l)^d\times Q_1)=
\mu_2(Q_1)=
\L^d(\psi_2^{-1}(Q_1))=
\L^d(B_1)  .  $$
Thus, we only have to find sets $A_{i,1}\subset A_i$ which satisfy the first formula of (1.8); since $\L^d$ is non-atomic and, by (1.7), 
$$\L^d(A_i)=\g(Q_i\times [-l,l)^d)\ge\g(Q_i\times Q_1)    $$
this is standard.  

Now, we find the sets $B_{i,1}$ which satisfy (1.9).  First of all we note that, by (1.8), 
$$\L^d(B_1\setminus\bigcup_{i\ge 2}A_{i,1})\ge
\L^d(A_{1,1})  .  $$
Since the $A_{i,1}$ are disjoint,we also have that 
$B_1\cap A_{1,1}$ does not intersect $A_{i,1}$ for $i\ge 2$; moreover, $\L^d(B_1\cap A_{1,1})\le\L^d(A_{1,1})$. Thus, we can find $B_{1,1}\subset B_1$ such that

\noindent $a$) $B_{1,1}\supset A_{1,1}\cap B_1$,

\noindent $b$) $\L^d(B_{1,1})=\L^d(A_{1,1})$,

\noindent $c$) $B_{1,1}$ is disjoint from 
$A_{i,1}$ for $i\ge 2$. 

Point $c$) follows by the last formula: in 
$B_1\setminus\bigcup_{i\ge 2}A_{i,1}$ there is enough space to accommodate a $B_{1,1}$ satisfying $b$). 

We show the next step of the induction, namely how to find 
$B_{2,1}$. By (1.8) and the aforesaid,
$$\L^d\left(
B_1\setminus\left( B_{1,1}\cup\bigcup_{i\not=2}A_{i,1}\right)
\right)\ge
\L^d(A_{2,1})  .  $$
Using this, we can find $B_{2,1}\subset B_1$ such that 

\noindent $a^\prime$) $B_{2,1}\supset A_{2,1}\cap B_1$,

\noindent $b^\prime$) $\L^d(B_{2,1})=\L^d(A_{2,1})$,

\noindent $c^\prime$) $B_{2,1}$ is disjoint from $B_{1,1}$ and from $A_{i,1}$ for 
$i\not=2$. 

Iterating, we get the sets $B_{i,1}$; the first, second, fourth and fifth formulas of (1.9) follow by construction, the third one by the first formula of (1.9), (1.8) and the fact that the $B_{i,1}$ are disjoint. 

\noindent{\bf Step 5.} In this step, we define $h_n$ on the first row of cubes: we want to find an invertible, bi-measurable map 
$\hat h_1$ which preserve Lebesgue measure and such that, for all $i$, 
$$\left\{
\matrix{
\hat h_1(x)=x\txt{if}x\not\in\bigcup_i(A_{i,1}\cup B_{i,1})\cr
(\psi_1\circ\hat h_1,\psi_2)^{-1}(Q_i\times Q_1)=B_{i,1} . 
}
\right.    \eqno (1.10)$$
Before proving this, note that $\L^d(B_{i,1})=\g(Q_i\times Q_1)$ by (1.8) and (1.9); this and (1.10) proves that (1.5) holds for the first row of cubes $\{ Q_i\times Q_1 \}_i$. The other rows will follow by induction, as we shall see in step 6. 

We prove (1.10). First of all, there are invertible maps 
$\fun{\phi_i}{B_{i,1}}{A_{i,1}}$ which preserve Lebesgue measure and which are the identity on $A_{i,1}\cap B_{i,1}$. This is easy to do: we set $\phi_i(x)=x$ on $A_{i,1}\cap B_{i,1}$; then, we use theorem 15.5.16 of [22] to get an invertible, measure-preserving map $\phi_i$ from $B_{i,1}\setminus A_{i,1}$ to 
$A_{i,1}\setminus B_{i,1}$; recall that these sets have the same Lebesgue measure by the first one of (1.9). 

Next, we glue together the maps $\phi_i$ in the following way:
$$\hat h_1(x)=\left\{
\eqalign{
x&\txt{if}x\not\in\bigcup_i(A_{i,1}\cup B_{i,1})\cr
\phi_i(x)&\txt{if}x\in B_{i,1}\cr
\phi_i^{-1}(x)&\txt{if}x\in A_{i,1}  . 
}
\right.    $$
The definition is well-posed: since by (1.9) the $B_{i,1}$ are disjoint, and since we saw above that the $A_{i,1}$ are disjoint, the only possible conflict is when $x\in B_{i,1}\cap A_{j,1}$. But then by (1.9) $j=i$; now on $B_{i,1}\cap A_{i,1}$ $\phi_i$ and 
$\phi_i^{-1}$ coincide, since both are the identity on this set. 

To check (1.10), we begin to note that its first formula comes straight from the definition of $\hat h_1$. As for the second one, if 
$x\in(\psi_1\circ\hat h_1,\psi_2)^{-1}(Q_i\times Q_1)$, then 
$x\in\psi_2^{-1}(Q_1)=B_1$ and 
$\hat h_1(x)\in\psi_1^{-1}(Q_i)=A_i$. Now $B_1$ is partitioned by the $B_{j,1}$ and the only $B_{j,1}$ which $\hat h_1$ sends to $A_i$ is $B_{i,1}$. Thus, $x\in B_{i,1}$, proving that 
$(\psi_1\circ\hat h_1,\psi_2)^{-1}(Q_i\times Q_1)= B_{i,1}$. 

\noindent{\bf Step 6.} We saw above that (1.5) follows if we show (1.10) for all the other rows; we do this by iteration. By the last step, the pre-image of $\cup_i(Q_i\times Q_1)$ by 
$(\psi_1\circ\hat h_1,\psi_2)$ is $B_1$. We want to adjust the second row of cubes without touching $B_1$. To do this, we restrict $(\psi_1\circ\hat h_1,\psi_2)$ to $B_1^c$; its image will fall in 
$$\bigcup_{j\not=1}(Q_i\times Q_j)   .  $$
Now we apply the procedure of the first step to the second row, i. e. to $\{ Q_i\times Q_2 \}_i$ and to $(\psi_1\circ\hat h_1,\psi_2)$. We get a map $\hat h_2$ from $B_1^c$ to itself such that 
$(\psi_1\circ\hat h_1\circ\hat h_2,\psi_2)$ satisfies (1.5) for $j=2$. Now we extend $\hat h_2$ to be the identity on $B_1$, and we get that $(\psi_1\circ\hat h_1\circ\hat h_2,\psi_2)$ satisfies (1.5) for $j=1$ too. To close, it suffices to call $h_n$ the last step of the iteration, the one in which all the rows are settled. 

\fin

We can look at $\W_2$ on $\Pc$ keeping track of the action of 
$\R^d$ on $\T^d$. Let us define 
$$\fun{\pi_{\T^d}}{\T^d\times\R^d}{\T^d}    $$
as the projection on the first coordinate, and let us set 
$$\fun{\a}{\T^d\times\R^d}{\T^d},\qquad
\fun{\a}{(x,v)}{x+v}  .  $$
Let $\tilde\mu_1,\tilde\mu_2\in\Pc$; we say that 
$\g\in{\cal P}_2(\T^d\times\R^d)$ belongs to  
$\Psi(\tilde\mu_1,\tilde\mu_2)$ if 
$(\pi_{\T^d})\sharp\g=\tilde\mu_1$ and 
$\a_\sharp\g=\tilde\mu_2$; we leave to the reader the simple proof that 
$$\W_2^2(\tilde\mu_1,\tilde\mu_2)=
\min_{\g\in\Psi(\tilde\mu_1,\tilde\mu_2)}
\int_{\T^d\times\R^d}|v|^2\dr\g(x,v)  .  \eqno (1.11)$$
We denote by $\Psi_o(\tilde\mu_1,\tilde\mu_2)$ the set of minimals. 

In the following, we shall denote by $L^2_\mu$ a space of $L^2$ functions for the measure $\mu$; we shall omit the $\mu$ when it is the Lebesgue measure. 

Let now $\fun{G}{\Pc}{\R}$ be a function; we say that $G$ is differentiable at $\tilde\mu\in\Pc$ if there is a vector field 
$\xi\in L^2_{\tilde\mu}(\T^d,\R^d)$ such that
$$\left\vert
G(\tilde\nu)-G(\tilde\mu)-\int_{\T^d\times\R^d}\inn{\xi(x)}{v}\dr\g(x,v)
\right\vert   =  o(\W_2(\tilde\mu,\tilde\nu))$$
for all $\tilde\nu\in\Pc$ and all $\g\in\Psi_o(\tilde\mu,\tilde\nu)$; we have denoted by 
$\inn{\cdot}{\cdot}$ the inner product in $\R^d$. 

Following [15], we say that $G$ is strongly differentiable at 
$\tilde\mu$ if there is $k>0$ such that 
$$\left\vert
G(\tilde\nu)-G(\tilde\mu)-\int_{\T^d\times\R^d}\inn{\xi(x)}{v}\dr\g(x,v)
\right\vert   \le
k\int_{\T^d\times\R^d}|v|^2\dr\g(x,v)$$
for all $\tilde\nu\in\Pc$ and all $\g\in\Psi(\tilde\mu,\tilde\nu)$. Note that we don't restrict the transfer plan $\g$ to be in 
$\Psi_o(\tilde\mu,\tilde\nu)$; it is immediate that strong differentiability implies differentiability. Of course, there are parallel definitions of differentiability and strong differentiability in $\Pcc$, which we forego to state. 

If $\fun{G}{\Pc}{\R}$, we can define
$$\fun{\bar G}{\Pcc}{\R},\qquad
\bar G(\mu)=G(\pi_\sharp\mu)  .  \eqno (1.12)$$

\lem{1.2} Let $\fun{G}{\Pc}{\R}$ be strongly differentiable at 
$\tilde\mu$ and let $\fun{\bar G}{\Pcc}{\R}$ be defined as in (1.12). Then, $\bar G$ is strongly differentiable at any 
$\mu\in\Pcc$ such that $\pi_\sharp\mu=\tilde\mu$. 

Conversely, if $\fun{\bar G}{\Pcc}{\R}$ quotients to a map 
$\fun{G}{\Pc}{\R}$ and is strongly differentiable at $\mu$, then 
$G$ is strongly differentiable at $\tilde\mu=\pi_\sharp\mu$.

\proof We begin with the direct statement. Let  
$\tilde\xi\in L^2(\T^d,{\tilde\mu})$ be the derivative of $G$ at 
$\tilde\mu$; we define $\fun{\xi}{\R^d}{\R^d}$ by 
$\xi(y)=\tilde\xi(\pi(y))$. We assert that $\xi\in L^2(\R^d,\mu)$; indeed, since $\pi_\sharp\mu=\tilde\mu$ we get the equality below, while the inequality comes from the fact that 
$\tilde\xi\in L^2_{\tilde\mu}$.
$$\int_{\R^d}|\xi(x)|^2\dr\mu(x)=
\int_{\T^d}|\tilde\xi(x)|^2\dr\tilde\mu(x)<+\infty . $$

We prove that $\xi$ is the derivative of $\bar G$ at $\mu$. Let 
$\nu\in\Pcc$ project on $\tilde\nu\in\Pc$ and let 
$\g\in\Psi(\mu,\nu)$; if we define 
$\tilde\g=(\pi\times id)_\sharp\g$ we see easily that 
$\tilde\g\in\Psi(\tilde\mu,\tilde\nu)$. We disintegrate $\g$ as 
$\mu\otimes\g_x$ and $\tilde\g$ as $\tilde\mu\otimes\tilde\g_q$, where $\g_x$ and $\tilde\g_q$ are measures on $\R^d$. An easy check shows that, if $f\in C(\T^d\times\R^d)$ with 
$\frac{f(x,v)}{1+|v|^2}$ bounded, then
$$\int_{\R^d}\dr\mu(x)\int_{\R^d}
f(x,y)\dr\g_x(y)=
\int_{\T^d}\dr\tilde\mu(q)\int_{\R^d}
f(q,y)\dr\tilde\g_q(y)  .  $$

The first equality below comes from (1.12) and the disintegration of $\g$; the second one comes from the definition of $\xi$ using the fact that $\tilde\mu=\pi_\sharp\mu$ and the formula above. The third equality comes from the disintegration of $\tilde\g$. The first inequality comes from the fact that $G$ is strongly differentiable, while the last equality is obvious. 
$$\left\vert
\bar G(\nu)-\bar G(\mu)-
\int_{\R^d\times\R^d}\inn{\xi(x)}{v}\dr\g(x,v)
\right\vert  =  $$
$$\left\vert
G(\tilde\nu)-G(\tilde\mu)-
\inn{\int_{\R^d}\xi(x)\dr\mu(x)}{\int_{\R^d}v\dr\g_x(v)}
\right\vert  =$$
$$\left\vert
G(\tilde\nu)-G(\tilde\mu)-
\inn{\int_{\T^d}\tilde\xi(q)\dr\tilde\mu(q)
}{\int_{\R^d}v\dr\tilde\g_q(v)}
\right\vert  =$$
$$\left\vert
G(\tilde\nu)-G(\tilde\mu)-
\int_{\T^d\times\R^d}\inn{\tilde\xi(q)}{v}\dr\tilde\g(q,v)
\right\vert  \le  $$
$$k\int_{\T^d\times\R^d}|v|^2\dr\tilde\g(x,v)=
k\int_{\R^d\times\R^d}|v|^2\dr\g(x,v)   .  $$
Since this is the definition of strong differentiability in $\Pcc$, we are done.

We prove the converse. 

\noindent{\bf Step 1.} Let $\tilde\mu,\tilde\nu\in\Pc$, let 
$\mu\in\Pcc$ be such that $\pi_\sharp\mu=\tilde\mu$ and let 
$\tilde\g\in\Psi(\tilde\mu,\tilde\nu)$. Recall that we have defined a map $\fun{\a}{(x,v)}{x+v}$.  We assert that we can find 
$\g\in{\cal P}_2(\R^d\times\R^d)$ and $\nu\in\Pcc$ such that 

\noindent $a$) the first marginal of $\g$ is $\mu$,

\noindent $b$) $(\pi\times id)_\sharp\g=\tilde\g$ and

\noindent $c$) $\a_\sharp\g=\nu$ and 
$\pi_\sharp\nu=\tilde\nu$; in particular, $\g\in\Psi(\mu,\nu)$.

To find $\g$, we disintegrate $\mu$ as 
$\mu=\b_q\otimes\tilde\mu$, with $\b_q$ a probability measure on the fiber $\{ q+\Z^d \}$; in other words, 
$\b_q(z)\ge 0$ and 
$$\sum_{z\in\Z^d}\b_q(z)=1  .  $$
Then, we can define $\g$ by
$$\int_{\R^d\times\R^d}f(x,v)\dr\g(x,v)=
\int_{\T^d\times\R^d}
\left[
\sum_{z\in\Z^d}\b_q(z)f(q+z,v)
\right]
\dr\tilde\g(q,v) $$
for all continuous functions $\fun{f}{\R^d\times\R^d}{\R}$ such that 
$\frac{f(x,v)}{1+|v|^2}$ is bounded. 
Setting $\nu=\a_\sharp\g$ we easily check that $\g$ and $\nu$  satisfy $a$), $b$) and $c$). 

\noindent{\bf Step 2.} Let $\xi$ be the derivative of 
$\bar G$ at $\mu$; we assert that $\xi=\tilde\xi\circ\pi$, where 
$\tilde\xi$ is a vector field on $\T^d$.  This is easy to see: for instance, taking a vector field $\eta$ supported in a small ball 
$B(x_0,r)$ of $\R^d$, considering 
$\g_{\e,z}=\mu\otimes(id+\e\eta(\cdot+z))_\sharp\L^d$ for 
$z\in\Z^d$, setting $\nu_{\e,z}=\a_\sharp\g_{\e,z}$ and noting that 
$\bar G(\nu_{\e,z})$, which quotients on $\Pc$, depends on $z$ only through $\mu(B(z_0,r))$. 

\noindent{\bf End of the proof.} The two steps above yield the first equality below, while the inequality comes from the fact that 
$\bar G$ is strongly differentiable at $\mu$. 
$$\left\vert
G(\tilde\nu)-G(\tilde\mu)-
\int_{\T^d\times\R^d}\inn{\tilde\xi(q)}{v}\dr\tilde\g(q,v)
\right\vert   =  $$
$$\left\vert
\bar G(\nu)-\bar G(\mu)-
\int_{\R^d\times\R^d}\inn{\xi(x)}{v}\dr\g(x,v)
\right\vert\le
k\int_{\R^d\times\R^d}|v|^2\dr\g(x,v)=
k\int_{\T^d\times\R^d}|v|^2\dr\tilde\g(q,v)   .  $$

\fin

We shall denote by $H$ the group of all bi-measurable maps 
$\fun{h}{[0,1)^d}{[0,1)^d}$ which preserve Lebesgue measure; we also set $L^2_\Z\colon= L^2([0,1)^d,\Z^d)$, which is a group under addition. 

Given $\fun{G}{\Pc}{\R}$, we can define a function
$$\fun{\hat G}{M}{\R},\qquad
\hat G(\psi)=G(\pi_\sharp\circ\psi_\sharp\L^d)  .  \eqno (1.13)$$
Clearly, the map $\hat G$ defined above is $H$ and $L^2_\Z$-equivariant, i. e.
$$\hat G(\psi\circ h+z)=\hat G(\psi)\qquad
\forall(\psi,h,z)\in M\times H\times L^2_\Z  .  \eqno (1.14)$$
Going in the opposite direction, if $\fun{\hat G}{M}{\R}$ is a continuous map such that (1.14) holds, we can define
$$\fun{\bar G}{\Pcc}{\R},\qquad
\bar G(\mu)=\hat G(\psi)   \eqno (1.15)$$
where $\psi\in M$ is such that $\psi_\sharp\L^p=\mu$. We prove that $\bar G$ is well-defined on $\Pcc$: actually, we are going to see that $\bar G$ quotients to a function $G$ on ${\cal P}(\T^d)$. Indeed, if $\psi_1,\psi_2\in M$ are such that $\pi_\sharp(\psi_i)_\sharp\L^p=\tilde\mu\in{\cal P}(\T^d)$ for $i=1,2$, then it is standard (lemma 6.4 of [11] or lemma 1.1 above) that there are $h_n\in H$  and $z_n\in L^2_\Z$ such that 
$$||\psi_1-\psi_2\circ h_n-z_n||_M\tends 0
\txt{as}n\tends+\infty  .  $$
The equality below comes from (1.14), while the limit comes from the formula above and the continuity of $\hat G$. 
$$\hat G(\psi_1)-\hat G(\psi_2)=
\hat G(\psi_1)-\hat G(\psi_2\circ h_n+z_n)\tends 0 . $$
This proves that $\hat G$ is well defined; as for the differentiability of $\hat G$, we recall theorems 6.2 and 6.5 of [11]. 

\prop{1.3} Let $\fun{\hat G}{M}{\R}$ be continuous and let it satisfy (1.14). Then, the following happens. 

\noindent 1) If $\hat G$ is differentiable at $\psi$, then $\hat G$ is differentiable at $\eta$ for all $\eta\in M$ such that 
$\eta_\sharp\L^d=\psi_\sharp\L^d$. Moreover, the law of 
$D\hat G(\psi)$ does not depend on the choice of $\eta$. 

\noindent 2) Let us suppose that $\hat G$ is of class $C^1$ and let 
$\mu\in\Pcc$. Then, there is 
$\xi\in L^2_\mu(\R^d,\R^d)$ such that, for all $\psi$ satisfying 
$\psi_\sharp\L^d=\mu$, we have
$$D\hat G(\psi)(x)=\xi\circ\psi(x)\txt{for $\L^p$ a. e. $x$.}$$

\noindent 3) Let $\hat G\in C^2(M,\R)$ with a bounded second derivative and let it satisfy (1.14); then, the function $\fun{\bar G}{\Pcc}{\R}$ defined by (1.15) is strongly differentiable. By lemma 1.2 this implies that its quotient $G$ on $\Pc$ is strongly differentiable. 

\proof Point 1) is theorem 6.2 of [11], point 2 theorem 6.5. We prove the easy consequence 3).

We want to show that $\bar G$ is strongly differentiable at any 
$\mu\in\Pcc$. Thus, let $\nu\in\Pcc$ and let $\psi,\eta\in M$ be such that $\psi_\sharp\L^p=\mu$, $\eta_\sharp\L^p=\nu$; let 
$\l\in\Psi(\mu,\nu)$ and let $\xi$ be as in point 2) above. Let 
$\fun{\b}{(x,v)}{(x,x+v)}$; since $\l\in\Psi(\mu,\nu)$ it is easy to check that $\g\colon=\b_\sharp\l$ belongs to $\G(\mu,\nu)$. By formula (1.3) of lemma 1.1 we can find $h_n\in H$ such that 
$$\int_{[0,1)^d}|\psi(x)-\eta\circ h_n(x)|^2\dr x\tends
\int_{\R^d\times\R^d}|q-q^\prime|^2\dr\g(q,q^\prime)$$
or equivalently, setting 
$\l_n\colon=(\psi,\eta\circ h_n-\psi)_\sharp\L^d$, 
$$\int_{\R^d\times\R^d}|v|^2\dr\l_n(x,v)\tends
\int_{\R^d\times\R^d}|v|^2\dr\l(x,v)  .  \eqno (1.16)$$
We assert that
$$\int_{\R^d\times\R^d}
\inn{\xi(x)}{v}\dr\l_n(x,v)\tends
\int_{\R^d\times\R^d}
\inn{\xi(x)}{v}\dr\l(x,v)  .  \eqno (1.17)$$
Indeed, if $\xi$ were continuous, this would follow from (1.3). In the general case, we can find a continuous vector field 
$\xi^\prime$ such that $||\xi-\xi^\prime||_{L^2_\mu}<\e$;  
the first inequalities in the two formulas below are H\"older  while the second ones come from (1.16). 
$$\left\vert
\int_{\R^d\times\R^d}
\inn{\xi-\xi^\prime}{v}\dr\l_n(x,v)
\right\vert  \le
||\xi-\xi^\prime||_{L^2_\mu}
\left[
\int_{\R^d\times\R^d}
|v|^2\dr\l_n(x,v)
\right]^\2   \le M ||\xi-\xi^\prime||_{L^2_\mu}\le M\e  , $$
$$\left\vert
\int_{\R^d\times\R^d}
\inn{\xi-\xi^\prime}{v}\dr\l(x,v)
\right\vert  \le
||\xi-\xi^\prime||_{L^2_\mu}
\left[
\int_{\R^d\times\R^d}
|v|^2\dr\l(x,v)
\right]^\2   \le M ||\xi-\xi^\prime||_{L^2_\mu}\le M\e  .  $$
These two formulas imply the second inequality below; the third one follows from (1.3) taking $n$ large enough. 
$$\left\vert
\int_{\R^d\times\R^d}
\inn{\xi(x)}{v}\dr(\l_n-\l)(x,v)
\right\vert\le$$
$$\left\vert
\int_{\R^d\times\R^d}
\inn{\xi-\xi^\prime}{v}\dr(\l_n-\l)
\right\vert  +
\left\vert
\int_{\R^d\times\R^d}
\inn{\xi^\prime}{v}\dr(\l_n-\l)
\right\vert   \le  $$
$$2\e M+
\left\vert
\int_{\R^d\times\R^d}
\inn{\xi^\prime}{v}\dr(\l_n-\l)
\right\vert   \le  2\e M+\e  .  $$
This proves (1.17). By (1.17), there is $\e_n\tends 0$ such that the first inequality below holds. The second one follows if we take $k$ to be the $\sup$ of $\2||D^2\hat G||$, which is finite by hypothesis. The last inequality follows from (1.16). 
$$\left\vert
\bar G(\nu)-\bar G(\mu)-
\int_{\R^d\times\R^d}\inn{\xi(x)}{v}\dr\l(x,v)
\right\vert\le$$
$$\left\vert
\hat G(\eta\circ h_n)-\hat G(\psi)-
\int_{[0,1)^d}\inn{\xi(\psi(x))}{\eta\circ h_n(x)-\psi(x)}\dx
\right\vert   +\e_n\le  $$
$$k\int_{[0,1)^d}|\eta\circ h_n(x)-\psi(x)|^2\dx+\e_n\le
k\int_{\T^d\times\R^d}|v|^2\dr\l(x,v)+2\e_n  .  $$
Letting $n\tends+\infty$, we recover the definition of strong differentiability at $\mu$. 

\fin

In the opposite direction, we have the following.

\lem{1.4}  Let $\fun{G}{\Pc}{\R}$ be a function and let 
$\fun{\hat G}{M}{\R}$ be defined as in (1.13). Let us suppose that $G$ is strongly differentiable at $\tilde\mu\in\Pc$, let $\mu\in\Pcc$ be a representative of $\tilde\mu$ and let $\psi\in M$ such that 
$\psi_\sharp\L^d=\mu$. Then, $\hat G$ is differentiable at 
$\psi\circ h+z$ for all $(h,z)\in H\times L^2_\Z$, and
$$D\hat G(u\circ h+z)=D\hat G(u)\circ h  .  \eqno (1.18)$$

\proof We define $\fun{\bar G}{\Pcc}{\R}$ as in (1.12); by lemma 1.2, $\bar G$ is strongly differentiable at any representative $\mu$ of $\tilde\mu$. 

Let $\xi$ be the derivative of $\bar G$ at $\mu$ and let 
$\psi\in M$ be such that $(\psi)_\sharp\L^p=\mu$. Let $\eta\in M$ and let us set $\nu=\eta_\sharp\L^p$. 
If we define $\l=(\psi,\eta-\psi)_\sharp\L^p$, we get the first equality below. Now $\l\in\Psi(\mu,\nu)$ and $G$ is strongly differentiable at $\mu$ with differential $\xi$; for some $k>0$ this implies the inequality below, while the last equality comes from the definitions of $\hat G$ and $\l$.
$$k\int_{[0,1)^d}|\psi(x)-\eta(x)|^2\dr x=
k\int_{\T^p\times\R^d}|v|^2\dr\l(x,v)\ge$$
$$\left\vert
\bar G(\nu)-\bar G(\mu)-
\int_{\R^d\times\R^d}\inn{\xi(q)}{v}\dr\l(q,v)
\right\vert   =  $$
$$\left\vert
\hat G(\eta)-\hat G(\psi)-
\int_{[0,1)^d}\inn{\xi\circ\psi(x)}{\eta(x)-\psi(x)}\dr x
\right\vert   .  $$
The last formula implies that $\hat G$ is differentiable at 
$\psi$. 

As for point 2), this is a general property of equivariant functions: if $T_h$ is a set of bounded linear operators from $M$ to $M$  having the group property
$$T_{h_1}\circ T_{h_2}=T_{h_1h_2}$$
then it is standard that 
$$D\hat G(T_h u)=[T_{h^{-1}}^T D\hat G(u)]  $$
where $A^T$ denotes the adjoint operator of $A$. Setting 
$T_h u\colon=u\circ h$ and substituting, we get (1.18). 

\fin

\vskip 2pc
\centerline{\bf \S 2}
\centerline{\bf Assumptions on the potential and the final condition}
\vskip 1pc

We recall the assumptions used in [15] from section 6 onward. 

We begin to suppose that we are given 
$U^0,U^1,\phi\in C^3(\T^d)$ such that the lifts of $\phi$ and $U^1$ to $\R^d$ are even.

Our potential is the function $\fun{\fc}{\Pc}{\R}$ defined by
$$\fc(\mu)=\2\int_{\T^d}(\phi\ast\mu)(z)\dr\mu(z)=
\2\int_{\T^d\times\T^d}\phi(z-z^\prime)\dr\mu(z)\dr\mu(z^\prime)  
$$
where the symbol $\ast$ denotes, as usual, convolution. The final condition is the function $\fun{\uc_0}{\Pc}{\R}$ given by
$$\uc_0(\mu)=
\int_{\T^d}[
U^0(z)+\2 (U^1\ast\mu)(z)
]  \dr\mu(z)=$$
$$\int_{\T^d\times\T^d}[
U^0(z)+\2 U^1(z-z^\prime)
]  \dr\mu(z)\dr\mu(z^\prime)  .  $$
It is shown in [15] that $\fc$ and $\uc$ are strongly differentiable. 

We recall from the introduction that we denote by $\dr$ the differential of functions on $\Pc$, by $D$ and $\nabla$ that of functions on $M$ and on $\R^d$ respectively. 

Always by [15], we have that
$$\dr\fc(\mu)=\nabla F(q,\mu)
\txt{and}
\dr\uc_0(\mu)=\nabla u_0(q,\mu)$$
where 
$$F(q,\mu)=(\phi\ast\mu)(q)
\txt{and}
u_0(q,\mu)=U^0(q)+(U^1\ast\mu)(q)  .   $$
By (1.13), $\fc$ and $\uc$ induce functions $\hat\fc$ and 
$\hat\uc_0$ on $M$; by the definition of push-forward we see that, if $\s\in M$,
$$\hat\fc(\s)=\2
\int_{[0,1)^d\times[0,1)^d}
\phi[\s(x)-\s(y)]\dr x\dr y , \eqno (2.1)_a$$
$$\hat\uc_0(\s)=
\int_{[0,1)^d\times[0,1)^d}\{
U^0(\s(x))+\2 U^1[
\s(x)-\s(y)
]
\}   \dr x\dr y   .    \eqno (2.1)_b$$
Also the functions $F$ and $u_0$ extend to parametrizations:
$$\fun{\hat F}{\R^d\times M}{\R^d},\qquad
\hat F(q,\s)=\int_{[0,1)^d}\phi[q-\s(x)]\dr x , \eqno (2.2)_a$$
$$\fun{\hat u_0}{\R^d\times M}{\R^d},\qquad
\hat u_0(q,\s)= U^0(q)+\int_{[0,1)^d} U^1[q-\s(x)]\dr x .   
\eqno (2.2)_b$$

We forego the proof of the following lemma, which follows from our hypotheses on $\phi$, $U^0$, $U^1$ and standard facts about the Nemitsky operators (see for instance [2]). 

\lem{2.1} Let $\fun{\hat\fc,\hat\uc_0}{M}{\R}$ be defined as in (2.1), let $\hat F,\hat u_0$ be as in (2.2). Then, 
$\hat\fc$ and $\hat\uc_0$ are functions of class $C^3$ on $M$. Denoting by  $\inn{\cdot}{\cdot}$ and by $\inn{\cdot}{\cdot}_M$ the inner products in $\R^d$ and in $M$ respectively, we have that
$$D\hat\fc(\s)\psi=
\int_{[0,1)^d\times[0,1)^d}
\inn{\nabla\phi[\s(x)-\s(y)]}{\psi(x)}\dr x\dr y=
\inn{\nabla\hat F(\s(\cdot),\s)}{\psi}_M$$
and
$$D\hat\uc_0(\s)\psi=
\int_{[0,1)^d\times[0,1)^d}
\inn{\nabla U^0(\s(x))+\nabla U^1[\s(x)-\s(y)]}{\psi(x)}\dr x\dr y=
\inn{\nabla\hat u_0(\s(\cdot),\s)}{\psi}_M . $$
In other words, $D\hat\fc(\s)$ is represented by the function 
$\nabla\hat F(\s(\cdot),\s)\in M$,   $D\hat\uc_0(\s)$ by the funtion $\nabla\hat u_0(\s(\cdot),\s)\in M$. 
The functions $\hat F$ and $\hat u_0$
are of class $C^3$ in both variables, with bounded first, second and third derivatives. Moreover, $\hat F$ and $\hat u_0$ are 
$\Z^d$-equivariant in the first variable; they are also $L^2_\Z$ and $H$-equivariant in the second one. 

\rm

\vskip 2pc

\centerline{\bf \S 3}
\centerline{\bf Minima on short time intervals}
\vskip 1pc

In lemmas 3.2-3.5 below, we recall the method of [10] for the minimals of the value function; in lemma 3.1, we prove that the value functions on measures and on parametrizations coincide. 

\vskip 1pc

\noindent{\bf Definitions.} Let $\fun{\mu}{(t,0)}{\Pc}$ be a curve of measures satisfying, in the weak sense (the precise definition is in the proof of lemma 3.1 below), the continuity equation
$$\partial_s\mu_s+\div(X\mu_s)=0  \eqno ( 3.1) $$
for a drift $X\in L^2((t,0)\times\T^d,\L^1\otimes\mu_t)$. We define the augmented action of $(\mu_s,X)$ as 
$$\ac(t,\mu_s,X)=
\int_t^0[
\2||X(s,\cdot)||_{L^2_{\mu_s}}^2-\fc(\mu_s)
]  \dr s   +  \uc_0(\mu_0)  .  $$
The value function on $\Pc$ is defined by
$$\fun{\uc}{(-\infty,0]\times\Pc}{\R},\qquad
\uc(t,\bar\mu)=\inf\ac(t,\mu_s,X)  \eqno (3.2)$$
where the $\inf$ is over all paths $(\mu_s,X)$ which satisfy (3.1) and such that $\mu_t=\bar\mu$. We are not going to need this, but the $\inf$ is actually a minimum. 

Augmented action and value function lift in a natural way to the space $M$. Given $t\le 0$ and  a curve $\s\in AC((t,0),M)$, we can define
$$\hat\ac(t,\s)=
\int_t^0[
\2||\dot\s_s||_M^2-\hat\fc(\s_s)
]  \dr s+\hat\uc_0(\s_0)  .  $$
For $t\le 0$ and $\psi\in M$, we set 
$$\hat\uc(t,\psi)=\inf\{
\hat\ac(t,\s)\st\s\in AC((t,0),M)\txt{and}\s_t=\psi
\}   .  $$

\lem{3.1} Let $\uc$ and $\hat\uc$ be defined as above. Then, the following holds. 

\noindent 1) The function $\hat\uc$ is  continuous. Moreover, it is $H$ and $L^2_\Z$-equivariant, i. e.
$$\hat\uc(t,\psi)=\hat\uc(t,\psi\circ h+z)\qquad
\forall(t,\psi,h,z)\in (-\infty,0]\times M\times H\times L^2_\Z  .  $$

\noindent 2) Let $\tilde\mu\in\Pc$ and let $\psi\in M$ be such that 
$(\pi\circ\psi)_\sharp\L^d=\tilde\mu$. Then,
$$\uc(t,\tilde\mu)=\hat\uc(t,\psi)  .  $$

\proof Point 1) is easy to dispatch, since continuity is standard; we follow [18] for equivariance. If $\s_s$ is an AC curve with 
$\s_t=\psi$, $h\in H$ and $z\in L^2_\Z$, then 
$\tilde\s_s=\s_s\circ h+z$ is AC and satisfies 
$\tilde\s_t=\psi\circ h+z$; moreover, since the Lagrangian and 
$\hat\uc_0$ are $L^2_\Z$ and $H$-equivariant, we see immediately that 
$${\cal A}(t,\s)={\cal A}(t,\tilde\s)  .  $$
Clearly, this implies that $\hat\uc(t,\psi\circ h+z)\le\hat\uc(t,\psi)$; the opposite inequality is similar. 

As for point 2), we begin to prove that 
$$\hat\uc(t,\psi)\le\uc(t,\tilde\mu)  .  \eqno (3.3)$$
We assert that this follows if we show that, for any curve 
$(\mu_s,X)$ satisfying (3.1) with $\mu_t=\tilde\mu$ we can find 
$\s\in AC([t,0],M)$ such that 

\noindent $i$)  
$(\pi\circ\s_t)_\sharp\L^d=(\pi\circ\psi)_\sharp\L^d=\tilde\mu$, 

\noindent $ii$) ${\cal A}(t,\mu_s,X)=\hat{\cal A}(t,\s)$. 

Indeed, we saw after formula (1.15) that $i$) together with point 1) of this lemma implies that $\hat{\uc}(t,\s_0)=\hat{\uc}(t,\psi)$; since $ii$) implies that $\hat\uc(t,\s_0)\le\uc(t,\tilde\mu)$ , formula (3.3) follows. 

Thus, let $(\mu_s,X)$ be a weak solution of (3.1) with 
$\mu_t=\tilde\mu$. By proposition 4.21 of [5] (or theorem 8.2.1 of [4]) there is a measure $\Xi$ on $C([t,0],\T^d)$ such that, denoting by 
$\fun{\eta_s}{C([t,0],\T^d)}{\T^d}$ the evaluation map 
$\fun{\eta_s}{\g}{\g_s}$, we have  
$$(\eta_s)_\sharp\Xi=\mu_s \txt{for all} s\in[t,0]  .  
\eqno (3.4)$$ 
Moreover, $\Xi$ concentrates on absolutely continuous curves and 
$$\int_{C([a,b],\T^d)}\dr\Xi(\g)
\int_t^0|\dot\g(s)|^2\dr s=\int_t^0||X(s,x)||_{L^2_{\mu_s}}^2\dr s  .  
\eqno (3.5)$$
It is standard (see for instance theorem 15.5.16 of [22]) that there is a Borel map $\fun{B}{[0,1)^d}{C([t,0],\T^d)}$ such that 
$\Xi=B_\sharp\L^d$. We set 
$$\s_s(x)=B(x)(s)=\eta_s\circ B(x)  .  $$
Now point $i$) follows from (3.4), since 
$(\s_t)_\sharp\L^d=(\eta_t\circ B)_\sharp\L^d=
(\eta_t)_\sharp\Xi=\mu_t$. We prove point $ii$). 

The first equality below is the definition of ${\cal A}$, the second one is implied by (3.4) and (3.5) while the third one follows because $\Xi=B_\sharp\L^d$ and  
$(\eta_0)_\sharp\Xi=\mu_0=(\s_0)_\sharp\L^d$. The last equality is the definition of $\hat{\cal A}$.
$${\cal A}(t,\mu_s,X)=\int_t^0\left[
\2||X(s,\cdot)||^2_{L^2_{\mu_s}}-
\2\int_{\T^d\times\T^d}\phi(q-q^\prime)
\dr\mu_s(q)\dr\mu_s(q^\prime)
\right]  \dr s+\uc_0(\mu_0)=$$
$$\int_t^0\dr s
\left[
\int_{C([a,b],\T^d)}
\2|\dot\g(s)|^2\dr\Xi(\g)-
\2\int_{C([a,b],\T^d)\times C([a,b],\T^d)}
\phi(\g(s)-\g^\prime(s))\dr\Xi(\g)\dr\Xi(\g^\prime)
\right]
+$$
$$+\uc_0((\eta_0)_\sharp\Xi)=
\int_t^0
\left[
\2||\dot\s_s||_M^2\dr s-
\int_t^0\hat\fc(\s_s)
\dr s
\right]     +\hat\uc_0(\s_0)=\hat{\cal A}(t,\s)  .  $$
To prove the inequality opposite to (3.3), we let 
$\s\in AC((t,0),M)$ with $\s_0=\psi$ and we define 
$$\mu_s=(\pi\circ\s_s)_\sharp\L^d\txt{for}s\in(t,0)  .  
\eqno (3.6)$$
We want to show 

\noindent a) that $\mu$ satisfies (3.1) for a suitable drift $X$ and

\noindent b) that the augmented action of $(\mu_s,X)$ isn't larger  than the augmented action of $\s$. 

Clearly, a) and b) imply the inequality opposite to (3.3), from which the thesis follows. We begin with a): the idea is that $X(s,q)$ is the average of the velocities $\dot\s_s(x)$ of the curves which satisfy $\s_s(x)=q$. 

The measure $\L^1\otimes(\pi\circ\s_s,\dot\s_s)_\sharp\L^d$ on 
$[t,0]\times\T^d\times\R^d$ has marginal 
$\L^1\otimes(\pi\circ\s_s)_\sharp\L^d$ on $[t,0]\times\T^d$; we disintegrate 
$\L^1\otimes(\pi\circ\s_s,\dot\s_s)_\sharp\L^d=
\L^1\otimes(\pi\circ\s_s)_\sharp\L^d\otimes\nu_{s,q}$ 
where $\nu_{s,q}$ is a measure on $\R^d$, depending in a Borel way on $(s,q)\in[t,0]\times\T^d$. In other words, if 
$f\in C(\T^d\times\R^d)$ is such that $\frac{|f(x,v)|}{1+|v|^2}$ is bounded, then the first equality below holds for $\L^1$ a. e. 
$s\in[a,b]$; the second equality comes from (3.6). 
$$\int_{[0,1)^d}f(\s_s(x),\dot\s_s(x))\dr x=
\int_{[0,1)^d}\dr x\int_{\R^d}f(\s_s(x),v)\dr\nu_{s,\s_s(x)}(v)=
\int_{\T^d}\dr\mu_s(q)\int_{\R^d}f(q,v)\dr\nu_{s,q}(v)
 .   \eqno (3.7)$$
We set 
$$X(s,q)=\int_{\R^d}v\dr\nu_{s,q}(v)$$
Let now $\phi\in C^\infty_c((t,0)\times\T^d)$; the first equality below comes from (3.6), the second one from the definition of $X$ and the third one from (3.7). The last equality follows since 
$\phi$ has compact support in $(t,0)\times\T^d$.
$$\int_t^0\dr s\int_{\T^d}[
\partial_s\phi(s,q)+\inn{\nabla\phi(s,q)}{X(s,q)}
]\dr\mu_s(q)=$$
$$\int_t^0\dr s
\int_{[0,1)^d}[
\partial_s\phi(s,\s_s(x))+
\inn{\nabla\phi(s,\s_s(x))}{X(s,\s_s(x))}
]  \dr x=$$
$$\int_t^0\dr s
\int_{[0,1)^d}\left[
\partial_s\phi(s,\s_s(x))+
\inn{\nabla\phi(s,\s_s(x))}{\int_{\R^d}v\dr\nu_{s,\s_s(x)}(v)}
\right]  \dr x=$$
$$\int_t^0\dr s\int_{[0,1)^d}[
\partial_s\phi(s,\s_s(x))+
\inn{\nabla\phi(s,\s_s(x))}{\dot\s_s(x)}
]  \dr x=$$
$$\int_t^0\left[
\frac{\dr}{\dr s}\int_{[0,1)^d}\phi(s,\s_s(x))\dr x
\right]  \dr s=0  .  $$
This means that $(\mu_s,X)$ is a weak solution of (3.1), i. e. point a) holds.  

As for b), it is the same calculation, up to the use of Jensen's inequality:
$$\int_t^0\left[
\2\int_{\T^d}|X(s,q)|^2\dr\mu_s(q)-\fc(\mu_s)
\right]    \dr s+\uc_0(\mu_0)\le$$
$$\int_t^0\left[
\2\int_{\R^d}|v|^2\dr\nu_{s,q}(v)-\hat\fc(\s_s)
\right]  \dr s+\hat\uc(\s_0)=
\int_t^0\left[
\2||\dot\s_s||_M^2-\hat\fc(\s_s)
\right]    \dr s+\hat\uc(\s_0)    .  $$

\fin

Secured by the last lemma, from now on we shall concentrate on 
$\hat\ac$ and $\hat\uc$.

\vskip 1pc

\noindent{\bf Definition.} By $H^1_M(t,0)$ we denote the space of the maps $\s\in AC((t,0),M)$ such that
$$||\s||^2_{H^1_M}\colon=
||\s_t||_M^2+\int_t^0||\dot\s_s||^2_M\dr s<+\infty  .  $$
It is standard ([1]) that this is a Hilbert space for the inner product
$$\inn{\s}{\eta}_{H^1_M}\colon=
\inn{\s_t}{\eta_t}_M+
\int_t^0\inn{\dot\s_s}{\dot\eta_s}\dr s  .  $$
We recall the Poincar\'e-Wirtinger inequality
$$\sup_{s\in(t,0)}||\s_s||_M\le
||\s_t||_M+|t|^\2\cdot||\s||_{H^1_M}  .  $$

\lem{3.2} For $t<0$, let us consider the functional
$$\fun{I}{H^1_M(t,0)}{\R},\qquad
\fun{I}{\s}{\hat\ac(t,\s)}    $$
where the augmented action $\hat\ac$ has been defined at the beginning of this section. 
Then, the following points hold.

1) The functional $I$ is of class $C^1$ on $H^1_M(t,0)$. For 
$\hat F$ and $\hat u_0$ defined as in (2.2), we have 
$$I^\prime(\s)(h)=
\int_t^0[
\inn{\dot\s_s}{\dot h_s}_M-
\inn{\nabla\hat F(\s_s(\cdot),\s_s)}{h_s}_M
]  \dr s+
\inn{\nabla\hat u(\s_0(\cdot),\s_0))}{h_0}_M  =$$
$$\int_t^0\inn{\dot\s_s}{\dot h_s}_M\dr s-
\int_t^0\dr s
\int_{[0,1)^d\times[0,1)^d}\inn{\nabla\phi(\s_s(x)-\s_s(y))}{h_s(x)}
\dr x\dr y+$$
$$\int_{[0,1)^d\times[0,1)^d}
\inn{\nabla U^0(\s_0(x))+\nabla U^1(\s_0(x)-\s_0(y))}{h_0(x)}
\dr x\dr y   .  
\eqno (3.8)$$
To explain the notation, we recall that $\nabla\hat F(\cdot,\s_s)$ is a $C^2$ function from $\T^d$ to $\R^d$ and thus 
$\nabla\hat F(\s_s(\cdot),\s_s)\in M$. 

\noindent 2) Let $\s\in H^1_M(t,0)$ be minimal in the definition of 
$\hat\uc(t,\psi)$; then, $\s$ solves 
$$\left\{
\eqalign{
\ddot\s_s(x)&=-(\nabla\phi\ast\mu_s)(\s_s(x))=
-\nabla\hat F(\s_s(x),\s_s)
\txt{for} s\in (t,0)\cr
\s_t(x)&=\psi(x)\cr
\dot\s_0(x)&=-\nabla U^0(\s_0(x))-
(\nabla U^1\ast\mu_0)(\s_0(x))=-\nabla\hat u_0(\s_0(x),\s_0)
}
\right.  \eqno (3.9)$$
where we have set $\mu_s=(\s_s)_\sharp\L^p$. The equalities are in the space $M$, i. e. they hold for a. e. $x\in[0,1)^d$. 

\proof Since the potential $\hat\fc$ and the final condition 
$\hat\uc$ are defined by (2.1), the proof of (3.8) is classical (see for instance [2]) and we forego it.

We recall the proof of point 2), which again is classical. Since $I$ is of class $C^1$ by point 1), if $\s$ minimizes $I$ under the constraint $\s_t=\psi$, then we must have that
$$I^\prime(\s)(h)=0\txt{for all} h\in H^1_M(t,0)
\txt{with}h_t=0  .  $$
Integrating by parts in (3.8), this implies that
$$\int_t^0
\inn{-\ddot\s_s-(\nabla\hat F(\s_s(\cdot),\s_s)}{h_s}_M\dr s+
\inn{\dot\s_0}{h_0}_M+
\inn{\nabla\hat u_0(\s_0(\cdot),\s_0)}{h_0}_M
=0  $$
for all $h\in H^1_M(t,0)$ with $h_t=0$. Clearly, this implies the first and third formulas of (3.9), while the second one comes from the boundary conditions on the minimal $\s$. 

\fin

Finding minima of $I$ is a delicate proposition (see for instance [21]) because Tonelli's theorem does not apply to the infinite-dimensional space $M$. However, in our case the implicit function theorem comes to the rescue: in the next three lemmas we recall the approach of [10] in our situation. In the next lemma, we denote by $B_X(\psi,r)$ the ball in $X$ of radius $r$ and centered in $\psi$. 

\lem{3.3} There are $T,r>0$ such that the following holds. Let 
$t\in[-T,0]$, and let $\psi\in M$; we shall denote by $\psi$ both the element of $M$ and the function of $H^1_M(t,0)$ constantly equal to $\psi$. 

\noindent 1) There is a unique function 
$\s^{(t,\psi)}\in C^1([-T,0],M)$ such that

\noindent $i$) $\s^{(t,\psi)}_s\in B_M(\psi,r)$ for $s\in[-T,0]$, and

\noindent $ii$) $\s^{(t,\psi)}$ satisfies (3.9).

By the Poincar\'e-Wirtinger inequality, this implies that 
(3.9) has a unique solution in $B_{H^1_M(-T,0)}(\psi,r^\prime)$ for some $r^\prime>0$.  

\noindent 2) The map
$$\fun{\Phi}{[-T,0]\times M}{H^1_M(-T,0)},\qquad
\fun{\Phi}{(t,\psi)}{\s^{(t,\psi)}}$$
is of class $C^2$ and equivariant, i. e. 
$\s^{(t,\psi\circ h+z)}=\si\circ h+z$ for all $h\in H$ and 
$z\in L^2_\Z$. 

\proof Let us consider the map
$$\fun{\Sigma}{[-T,0]\times M}{M},\qquad
\fun{\Sigma}{(s,\tilde\psi)}{\s_s}$$
where $\s_s$ solves the Cauchy problem
$$\left\{
\eqalign{
\ddot\s_s(x)&=-\nabla\hat F(\s_s(x),\s_s)\cr
\s_0&=\tilde\psi\cr
\dot\s_0(x)&=-\nabla\hat u_0(\s_0(x),\s_0)=
-\nabla\hat u_0(\tilde\psi(x),\psi)
}
\right.   \eqno (3.10)$$
for the functions $\hat F$ and $\hat u$ which have been defined in (2.2). Since these two functions are of class $C^3$ by lemma 2.1, their gradients are in $C^2$ and the map $\Sigma$ is of class $C^2$ by the continuous dependence theorem. 

\noindent{\bf Step 1.} We assert that points 1) and 2) follow if we show that there is a $C^2$ function 
$\fun{\tilde\psi}{[-T,0]\times M}{M}$ which is, for all $\psi\in M$, the unique solution in $B(\psi,r)$ of 
$$\Sigma(t,\tilde\psi(t,\psi))=\psi  .  \eqno (3.11)$$
Indeed, if this holds we can set 
$$\s_s^{(t,\psi)}=\Sigma(s,\tilde\psi(t,\psi))  \eqno (3.12)$$
and (3.11) immediately implies that 
$$\s_t^{(t,\psi)}=\psi$$
i. e. $\s^{(t,\psi)}$ satisfies the second equation of (3.9).

Moreover, the map $\fun{}{(t,\psi,s)}{\s_s^{(t,\psi)}}$ is of class $C^2$ because of (3.12) and the fact that $\Sigma$ and 
$\tilde\psi$ are of class $C^2$; in particular, 
$\s^{(t,\psi)}\in H^1_M(-T,0)$. The map $\si$ solves the first equation of (3.9) because 
$\fun{}{s}{\Sigma(s,\tilde\psi(t,\psi))}$ solves it by the definition of 
$\Sigma$. Finally, $\si$ satisfies the third equation of (3.9) simply because it satisfies the third equation of (3.10). Uniqueness follows because, if (3.9) had two different solutions in 
$B_M(\psi,r)$, then also (3.11) would have two different solutions in $B_M(\psi,r)$, and we are supposing that this is not the case. 

We prove the last assertion of the lemma, equivariance. Recall that $\hat F$ and $\hat u_0$ are $H$ and $L^2_\Z$-equivariant; in particular, if $\si$ satisfies (3.9) and $(h,z)\in H\times L^2_\Z$, then also $\si\circ h+z$ satisfies (3.9) for the initial condition 
$\psi\circ h+z$. By the uniqueness of point 1), this implies that 
$\s^{(t,\psi\circ h+z)}=\si+z$ for all $h\in H$ and $z\in L^2_\Z$. 

\noindent{\bf Step 2.} In this step and in the following ones, we check that we can apply the implicit function theorem to solve for 
$\psi$ in (3.11). 

First of all, we saw above that the map $\Sigma$ is 
$C^2$. By definition, $\Sigma(0,\psi)=\psi$ for all $\psi\in M$, which implies that 
$$D\Sigma(0,\psi_0)=Id\qquad
\forall\psi_0\in M  .  $$

Thus, the implicit function theorem yields the existence of a $C^2$ function $\tilde\psi(t,\psi)$ defined in 
$[-T_0,0]\times B_M(\psi_0,r)$ which solves (3.1). 

In step 3 below, we shall see that $T_0$ and $r$ do not depend on $\psi_0$; in step 4, we shall use the monodromy theorem to glue the local solutions into a solution defined globally on 
$[-T_0,0]\times M$. 

\noindent {\bf Step 3.} We prove that we can choose $T_0$ and 
$r$ independent on $\psi_0$. 

If we look at the proof of the implicit function theorem, we see that $T_0,r>0$ must be chosen in order that the Lipschitz constant of 
$\fun{}{\psi}{\Sigma(t,\psi)-\psi}$ is smaller than, say, $\2$ in 
$[-T_0,0]\times B(\psi_0,r)$; by the Lagrange theorem, this follows if $|| D\Sigma(t,\psi)-Id ||\le\2$ in $[-T_0,0]\times B(\psi_0,r)$. This follows by a Taylor development, since we saw above that 
$D\Sigma(0,\psi) - Id=0$ for all $\psi$ and that 
$|| \partial_tD\Sigma(t,\psi) ||$ is bounded in $[-1,0]\times M$. 

\noindent{\bf Step 4.} By the last step, in each neighbourhood 
$[-T_0,0]\times B(\psi_0,r)$ we can define a function $\tilde\psi$ which satisfies (3.12); since $M$ is simply connected, we can use the monodromy theorem (see for instance theorem 1.8 of chapter 3 of [2]) to define globally a function 
$\fun{\tilde\psi}{[-T_0,0]\times M}{M}$ satisfying (3.11). 

\fin

\noindent{\bf Definition.} From now on, $\si_s$ will be defined as in the last lemma.

\vskip 1pc

Since the map $\fun{}{(t,\psi)}{\s^{(t,\psi)}}$ is of class $C^2$, the next lemma reduces to a classical computation ([10]) which we are only going to sketch; we continue in our practice of denoting by $D$ the derivative in the $M$ variable. 

\lem{3.4} We set
$$\hat\vc(t,\psi)=\int_t^0[
\2 ||\dot\si_s||_M^2-\hat\fc(\si_s)
]  \dr s+
\hat\uc_0(\si_0)  .  \eqno (3.13)$$
Then, $\hat\vc\in C^2([-T,0]\times M)$ and we have  
$$\left\{
\eqalign{
-\partial_t\hat\vc(t,\psi)+\2||D\hat\vc(t,\psi)||_M^2+\hat\fc(\psi)
&=0\txt{for}(t,\psi)\in[-T,0]\times M\cr
\hat\vc(0,\psi)&=\hat\uc_0(\psi)  . 
}
\right.  \eqno (3.14) $$
Moreover,
$$\dot\si_s=-D\hat\vc(s,\s^{(t,\psi)}_s)
\txt{for all}s,t\in[-T,0]  .  \eqno (3.15)$$

\proof  First of all, $\hat\vc\in C^2([-T,0]\times M)$ by point 2) of lemma 3.3. Next, we differentiate with respect to $\psi$ both terms of (3.13); after using (3.8) and (3.9) we get that
$$\dot\s_t^{(t,\psi)}=-D\hat\vc(t,\s_t^{(t,\psi)})=
-D\hat\vc(t,\psi) . \eqno (3.16)$$
Now we differentiate in (3.13) with respect to $t$; after an integration by parts, we get that 
$$\partial_t\hat\vc(t,\psi)=
-\2||\dot\s_t^{(t,\psi)}||^2_M+
\hat\fc(\si_t)+$$
$$\int_t^0\inn{
-\ddot\s_s^{(t,\psi)}-D\hat\fc(\s_s^{(t,\psi)})
}{
\partial_t\s_t^{(s,\psi)}
}_M   \dr s    +$$
$$\inn{
\dot\s_s^{(t,\psi)}
}{
\partial_t\s_s^{(t,\psi)}
}_M |_{s=t}^{s=0}+\inn{
D\hat\uc(\s_0^{(t,\psi)})
}{
\partial_t\s_0^{(t,\psi)}
}_M   .  $$
We note that the integral term is zero by the first equation of (3.9). Since $\si_t=\psi$ for all $t$, differentiating we get that
$$\partial_t\si_s|_{s=t}=-\dot\si_t  .   $$
Together with the last equation of (3.9), the last two equations imply that 
$$\partial_t\hat\vc(t,\psi)=
\2||\dot\s_t^{(t,\psi)}||^2_M+\hat\fc(\s_t^{(t,\psi)})  .  $$
Bt (3.16), this implies (3.14). 

Next, we assert that (3.15) follows from (3.16) if we show that, for all $t,s,\tau\in[-T,0]$, we have that 
$$\s^{(t,\psi)}_\tau=\s_\tau^{(s,\s_s^{(t,\psi)})} . \eqno (3.17)$$
To show the assertion, we denote by the dot the derivative in the 
$\tau$ variable; now (3.17) implies the first equality below, (3.16) the second one. 
$$\dot\s_\tau^{(t,\psi)}|_{\tau=s}=
\dot\s_\tau^{(s,\s_s^{(t,\psi)})}|_{\tau=s}=
-D\hat\vc(s,\si_s)  .  $$
To show (3.17), by the uniqueness of lemma 3.3 it suffices to show that $\fun{}{\tau}{\s^{(t,\psi)}_\tau}$ satisfies 
$$\left\{
\eqalign{
\ddot\s_\tau^{(t,\psi)}(x)&=
-\nabla\hat F(\s_\tau^{(t,\psi)}(x),\s_\tau^{(t,\psi)})\cr
\s_s^{(t,\psi)}(x)&=\s_s^{(t,\psi)}(x)\cr
\dot\s_0^{(t,\psi)}(x)&=
-\nabla\hat u_0(\s_0^{(t,\psi)}(x),\s_0^{(t,\psi)})
}
\right.$$
which is obvious since $\s^{(t,\psi)}$ satisfies (3.9).

\fin

\lem{3.5} Let $t\in[-T,0]$ and let $\psi\in M$. Then, 

\noindent 1) for all $s\in[-T,0]$,  $\si$ is the unique minimal in the definition of $\hat\uc(s,\si_s)$.

\noindent 2) $\hat\uc(t,\psi)=\hat\vc(t,\psi)$ for 
$(t,\psi)\in[-T,0]\times M$. 

\proof Point 2) follows immediately from point 1) and the definitions of $\hat\uc$ and $\hat\vc$; we recall the classical proof of [10] for point 1).   Let $\hat\vc$ be as in the last lemma and let us consider the functional
$$\fun{J_s}{H^1_M(t,0)}{\R},$$
$$\fun{J_s}{\s}{
\int_s^0[
\2||\dot\s_\tau||^2_M-\fc(\s_\tau)+\partial_\tau\hat\vc(\tau,\s_\tau)+
\inn{D\hat\vc(\tau,\s_\tau)}{\dot\s_\tau}_M
]  \dr  \tau
}      .   \eqno (3.18)$$
Since $\hat\vc$ is of class $C^2$ by lemma 3.4, we get the first equality below, while the second one follows from the second formula of (3.14) and the definition of $\hat{\cal A}$ at the beginning of this section. 
$$J_s(\s)=\int_s^0[
\2||\dot\s_\tau||_M^2-\fc(\s_\tau)
]\dr \tau +\hat\vc(0,\s_0)-\hat\vc(s,\s_s)=$$
$$\hat{\cal A}(s,\s)-\hat\vc(s,\s_s)  .  \eqno (3.19)$$
Thus, if we restrict to the curves $\s\in H^1_M(s,0)$ with 
$\s_s=\si_s$, minimizing $J_s$ is the same as minimizing 
$\hat{\cal A}(s,\s)$: the thesis follows if we check that $\si$ is minimal for $J_s$. Actually, we are going to show that the integrand of $J_s$ is constantly equal to its minimum along 
$(\tau,\si_\tau,\dot\si_\tau)$. 

Clearly, for all $(\tau,\eta)\in [-T,0]\times M$ the minimum of the Lagrangian of $J_s$ 
$$\fun{B_{\tau,\eta}}{M}{\R}$$
$$\fun{B_{\tau,\eta}}{\dot\l}{
\2||\dot\l||_M^2-\fc(\eta)+\partial_\tau\hat\vc(\tau,\eta)+
\inn{D_\eta\hat\vc(\tau,\eta)}{\dot\l}_M
}  $$
is attained at $\dot\l=-D_\eta\hat\vc(\tau,\eta)$; substituting this value into the expression for $B_{\tau,\eta}$ we get the inequality below, while the equality is the first formula of (3.14). 
$$B_{\tau,\eta}(\dot\l)\ge
-\2||D_\eta\hat\vc(\tau,\eta)||_M^2-
\fc(\eta)+\partial_\tau\hat\vc(\tau,\eta)=0 
\qquad\forall\dot\l\in M   . \eqno (3.20)$$
On the other side, (3.15) implies the second equality below, (3.14) the third one.
$$B_{\tau,\dot\si_\tau}(\dot\si_\tau)=
\2||\dot\si_\tau||_M^2-\fc(\si_\tau)+\partial_\tau\hat\vc(\tau,\si_\tau)+
\inn{D\hat\vc(\tau,\si_\tau)}{\dot\si_\tau}_M=$$
$$-\2||D\hat\vc(\tau,\si_\tau)||_M^2-\fc(\si_\tau)+
\partial_\tau\hat\vc(\tau,\si_\tau)=0  .  $$
The last two formulas imply that $\fun{}{\tau}{\si_\tau}$ minimizes $J_s$, as we wanted.

We prove uniqueness: by the aforesaid, if $\s_\tau$ minimizes, then the integrand of $J_s$ must be zero along $\s_\tau$. By (3.20), this implies that $\dot\s_\tau=-D\vc(\tau,\s_\tau)$. By (3.15) this implies that $\s_\tau$ and $\si_\tau$ satisfy the same differential equation; we recall from lemma 3.4 that 
$-D\hat\vc(t,\psi)$ is Lipschitz. Since $\s_s=\si_s$ by hypothesis, we get that $\s_\tau=\si_\tau$ for $\tau\in[-T,0]$ by the existence and uniqueness theorem. 

\fin

\vskip 2pc
\centerline{\bf\S 4}
\centerline{\bf The master equation}
\vskip 1pc

In this section, we are going to define the value function for the single particle; we shall see that it determines the movement of the whole pack and that it satisfies the master equation. 

\vskip 1pc

\noindent{\bf Definition.} We define
$$\fun{v}{[-T,0]\times\T^d\times[-T,0]\times M}{\R}  ,  $$
$$v(s,q|t,\psi)=
\min\left\{
\int_s^0[
\2|\dot y(\tau)|^2-\hat F(y(\tau),\si_\tau)
]\dr\tau +\hat u_0(y(0),\si_0)
\right\}  \eqno (4.1)$$
where the minimum (whose existence is guaranteed by Tonelli's theorem) is over all $y\in AC((s,0),\T^p)$ such that $y(s)=q$. In the notation for $v$ we have inaugurated the practice of placing the "parameters", in this case $(t,\psi)$, after the vertical slash. In other words, we are interested in the equation solved by $v$ in the first two variables. If we freeze $(t,\psi)$, then $v(s,q|t,\psi)$ is the value function of the particle $q$, given that the whole pack moves like $\si$. Thus, $v$ solves, in its first two variables, the Hamilton-Jacobi equation. 

\lem{4.1} Up to reducing $T$, the following holds. 

\noindent 1) For $s,t\in[-T,0]$, the minimum in the definition of 
$v(s,q|t,\psi)$ is attained on a unique function
$$\fun{}{\tau}{y(\tau|s,q,t,\psi)} . $$
Again, the parameters of the orbit (i. e. the initial conditions of the single particle and of the whole pack) are on the right of the vertical slash. 

\noindent 2) The map 
$$\fun{}{(\tau,s,q,t,\psi)}{y(\tau|s,q,t,\psi)}$$
is of class $C^2$.

\noindent 3) The value function
$$\fun{}{(s,q,t,\psi)}{v(s,q|t,\psi)}$$
is of class $C^2$ with bounded first and second derivatives. It is 
$\Z^d$-equivariant in the second variable, $H$ and $L^2_\Z$-equivariant in the fourth one. For all 
$(t,\psi)\in[-T,0]\times M$ it satisfies the Hamilton-Jacobi equation with time reversed
$$\left\{
\eqalign{
-\partial_sv(s,q|t,\psi)+\2|\nabla v(s,q|t,\psi)|^2+
\hat F(q,\si_s)&=0\quad (s,q)\in[-T,0]\times\T^d\cr
v(0,q|t,\psi)&=\hat u_0(q,\si_0)
}
\right.   \eqno (4.2)$$
in the classical sense. Recall that we denote the gradient in the 
$\T^p$ variable by $\nabla$, in the $M$ variable by $D$. 

\noindent 4) We have that, for $\L^p$ a. e. $x\in[0,1)^d$ and all 
$t,s,\tau\in[-T,0]$,
$$\dot y(\tau|s,\si_s(x),t,\psi)=\dot\si_\tau(x)=
-\nabla v(\tau,y(\tau|s,\si_s(x),t,\psi)|t,\psi)=
-D\hat\vc(\tau,\si_\tau)(x)  .  $$

\noindent 5) Let us define the function $S$ as the flow of 
$-\nabla v$, i. e. as 
$$S(s,q,\tau|t,\psi)=y(\tau)$$
where $y$ solves 
$$\left\{
\eqalign{
\dot y(\tau)&=-\nabla v(\tau,y(\tau)|t,\psi)\cr
y(s)&=q  .
}
\right.  \eqno (4.3)$$
Then, up to reducing $T$,  there is $D_2>0$ independent of 
$(s,q,\tau,t,\psi)\in[-T,0]\times\T^d\times[-T,0]^2\times M$ such that 
$$\frac{1}{D_2}\le
{\rm det}\frac{\partial S(s,q,\tau|t,\psi)}{\partial q}\le D_2  .  $$

\proof We fix $(t,\psi)$ as the initial condition of the whole pack; we consider the time dependent Lagrangian
$$\L(s,q,\dot q)=\2|\dot q|^2-\hat F(q,\si_s)$$
and the final condition
$$\fun{}{q}{\hat u_0(q,\si_0)}  .  $$
Note that, by lemma 2.1, $\L$ is $C^3$ in $(s,q,\dot q)$; it depends in a $C^2$ way on the parameters $(t,\psi)$ by lemma 3.3. Analogously, $\hat u_0$ is $C^3$ in the variable $q$ and 
$C^2$ in $(t,\psi)$. Now points 1), 2) and 3) follow by the argument of [10], which we have seen in lemmas 3.3, 3.4 and 3.5 above. 

As for point 4), formula (3.15) gives that, for all $\tau\in[-T,0]$, 
$$\dot\si_\tau(x)=-D\hat\vc(\tau,\si_\tau)(x)
\txt{for $\L^p$ a. e.} x\in [0,1)^d  .  $$
On the other side, with exactly the same proof we used for formula (3.15) we see that 
$$\dot y(\tau|s,\si_s(x),t,\psi)=
-\nabla v(\tau,y(\tau|s,\si_s(x),t,\psi)|t,\psi)
\txt{for}t,s,\tau\in[-T,0]  .  $$
Thus, it suffices to show the first equality of point 4). Classical Hamilton-Jacobi theory (which we recalled above in lemmas 3.3 to 3.5) implies that the minimizer
$$\fun{}{\tau}{y(\tau|s,q,t,\psi)}$$
satisfies
$$\left\{
\eqalign{
\frac{\dr^2}{\dr\tau^2}y(\tau|s,q,t,\psi)&=
-\nabla\hat F(y(\tau|s,q,t,\psi),\si_\tau)\cr
y(s|s,q,t,\psi)&=q\cr
\dot y(0|s,q,t,\psi)&=-\nabla\hat u_0(y(0|s,q,t,\psi),\si_0)  .
}
\right.  $$
If $q=\si_s(x)$ then, by (3.9), this is the same equation that is satisfied by $\fun{}{\tau}{\si_\tau(x)}$ for 
$\L^d$ a. e. $x\in[0,1)^d$; by the uniqueness of lemma 3.3 this implies the first equality of point 4). 

We prove point 5). Since $S(s,q,s|t,\psi)=q$ by definition, we see that $\partial_q S(s,q,s|t,\psi)=Id$; thus, point 5) follows if we show that the map 
$\fun{}{\tau}{\partial_q S(s,q,\tau|t,\psi)}$ is Lipschitz uniformly in 
$(s,q,\tau,t,\psi)$; in other words, we have to show that the norm of $\partial^2_{q\tau}S(s,q,\tau|t,\psi)$ is bounded. This follows easily by (4.3), the differentiable dependence theorem and point 3) of this lemma, which implies  
$$|\partial^2_{q,q}v(s,q|t,\psi)|\le M\qquad
\forall(s,q,t,\psi)\in[-T,0]\times\T^d\times[-T,0]\times M  .  $$

\fin

We can apply to the value function $v(s,q|t,\psi)$ a change of coordinates: namely, instead of seeing it as a function of 
$\si_t=\psi$, we can see it as a function of $\si_s$. In other words, we can define a function $u$ as 
$$u(s,q|\si_s)\colon=v(s,q|t,\psi)  .  $$
Equivalently, by (3.17) we get that, for $\psi\in M$, 
$\psi=\s_t^{(s,\s_s^{(t,\psi)})}$; setting $\eta=\s_s^{(t,\psi)}$ and substituting in the formula above, we get that
$$u(s,q|\eta)=v(s,q|t,\s_t^{(s,\eta)}) 
\txt{for all}t\in[-T,0],\quad\eta\in M   \eqno (4.4)$$
which incidentally proves that the definition of $u$ is well posed.
The first equality below comes from (4.4), since 
$\s_s^{(s,\psi)}=\psi$; the second one is (4.1).  
$$u(s,q|\psi)=v(s,q|s,\psi)=$$
$$\min\Big\{
\int_s^0[
\2|\dot y(\tau)|^2-\hat F(y(\tau),\s_\tau^{(s,\psi)})
]\dr\tau +\hat u_0(y(0),\s_0^{(s,\psi)})\st
y\in AC((s,0),\T^p),\quad y(s)=q
\Big\}  .   \eqno (4.5)$$

\lem{4.2} Let
$$\fun{u}{[-T,0]\times\T^d\times M}{\R}$$
be defined as in (4.4) or as in (4.5), which is the same. Then, 
$u$ is of class $C^2$ in all its variables and satisfies the master equation
$$-\partial_tu(t,q|\psi)+\2|\nabla u(t,q|\psi)|^2+
F(q,\psi)+
\inn{\nabla u(t,\psi(\cdot)|\psi)}{D u(t,q|\psi)}_M  =0  .  $$

\proof By (4.4), lemma 4.1 and the chain rule we get that $u$ is of class $C^2$ in all its variables. Since $\si_t=\psi$ for all $t$, differentiating we get that
$$\frac{\partial}{\partial s}\s_t^{(s,\psi)}|_{s=t}=
-\dot\si_t  .  \eqno (4.6)$$
The first equality of (4.5) implies the equalities below. 
$$D u(t,q|\psi)=
D v(t,q|t,\psi),\qquad
\nabla u(t,q|\psi)=
\nabla v(t,q|t,\psi).  \eqno (4.7)$$
The first equality below is point 4) of lemma 4.1, the second one comes from (4.7). 
$$\dot\si_t(x)=-\nabla v(t,\psi(x)|t,\psi)=-\nabla u(t,\psi(x)|\psi)   .    
\eqno (4.8)  $$
If we differentiate (4.4) in $s$, we get the first equality below; the second one comes from (4.2) and (4.6); the last one comes from (4.7) and (4.8).
$$\partial_su(s,q|\psi)|_{s=t}=
\partial_sv(s,q|t,\s_t^{(s,\psi)})|_{s=t}+
\inn{D v(s,q|t,\s_t^{(s,\psi)})}{\frac{\partial}{\partial s}\s_t^{(s,\psi)}}_M
|_{s=t}=$$
$$\2|\nabla v(t,q|t,\psi)|^2+\hat F(q,\psi)-
\inn{D v(t,q|t,\psi)}{\dot\si_t}_M=$$
$$\2|\nabla u(t,q|\psi)|^2+\hat F(q,\psi)+
\inn{D u(t,q|\psi)}{\nabla u(t,\psi(\cdot)|\psi)}_M  .  $$

\fin

\noindent{\bf End of the proof of theorem 1.} Point 1) follows from lemma 3.5; point 2) is point 2) of lemma 3.3; point 3) is lemma 4.2; point 4) follows from point 5) of lemma 4.1; point 5) is point 4) of lemma 4.1 and (4.7).

\fin

\noindent {\bf Remark.} By the results of section 1, $u(t,q|\psi)$ quotients to a function on measures which is strongly differentiable, with continuous derivative; it satisfies the master equation in the classical sense, i. e. taking derivatives at their face value.

\vskip 2pc
\centerline{\bf Bibliography}

\noindent [1] R. Adams, J. J. F. Fournier, Sobolev spaces, Academic Press, Singapore, 2009. 

\noindent [2] A. Ambrosetti, G. Prodi, A primer of nonlinear analysis, Cambridge University Press, Cambridge, 1995.

%\noindent [1] L. Ambrosio, W. Gangbo, Hamiltonian ODE's in the Wasserstein space of probability measures, Communications on Pure and Applied Math., {\bf 61}, 18-53, 2008.

\noindent [3] L. Ambrosio, J. Feng, On a class of first order Hamilton-Jacobi equations in metric space, preprint. 

\noindent [4] L. Ambrosio, N. Gigli, G. Savar\'e, Gradient Flows, Birkhaeuser, Basel, 2005.

\noindent [5] L. Ambrosio, N. Gigli, G. Savar\'e, Heat flow and calculus on metric measure spaces with Ricci curvature bounded below - the compact case. Analysis and numerics of Partial Differential Equations, 63-115, Springer, Milano, 2013.  

\noindent [6] R. J. Aumann, Markets with a continuum of traders, Econometrica, {\bf 32}, 39-50, 1964.

%\noindent[4] L. Ambrosio, G. Savar\'e, L. Zambotti, Existence and stability for Fokker-Planck equations with log-concave reference measures, preprint.

%\noindent [BB] P. Bernard, B. Buffoni, Optimal mass transportation and Mather theory, J. Eur. Math. Soc., {\bf 9}, 85-121, 2007.

%\noindent [10] U. Bessi, A time step approximation scheme for a viscous version of the Vlasov equation, to appear on Advances in Mathematics.

%\noindent [11] V. Bogachev, G. Da Prato, M. R\"ockner, Uniqueness for solutions of Fokker-Planck equations on infinite-dimensional spaces, J. Evol. Equ., {\bf 10}, 487-509, 2010.

%\noindent [15] G. Da Prato, F. Flandoli, M. R\"ockner, Uniqueness for continuity equations in Hilbert spaces with weakly differentiable drift, preprint.  

%\noindent [16] G. Da Prato, Introduction to stochastic analysis and Malliavin calculus, Pisa, 2007.

\noindent [7] A. Bensoussan, J. Frehse, P. Yam, The Master Equation in Mean Field Theory, J. Math. Pures Appl., {\bf 103}, 1441-1474, 2015. 

\noindent [8] A. Bensoussan, J. Frehse, P. Yam, On the interpretation of the Master Equation, Arxiv:1503.07754.

\noindent [9] R. Buckdahn, J. Li, S. Peng, C. Rainer, Mean Field Stochastic Differential Equations and associated PDE's, Arxiv:1407:1215.

\noindent [10] C. Caratheodory, Calculus of variations and partial differential equations of the first order, Chelsea, N. Y., 1989.

\noindent [11] P. Cardaliaguet, Notes on mean field games, from P. L. Lions' lectures at Coll\`ege de France, 

\noindent mimeographed notes.

\noindent [12] R. Carmona, F. Delarue, The Master equation for large population equilibria, Stochastic analysis and applications, Springer proceedings in Mathematical Statistics, {\bf 100}, 77-128, 2014. 

\noindent [13] J-F Chassagneux, D. Crisan, F. Delarue, A probabilistic approach to classical solutions of the Master Equation for large population equilibria, ArXiv:1411.3009v2.

%\noindent [D] R. L. Dobrushin, Vlasov equations, Functional analysis and its applications, {\bf 13}, 45-58, 1979.

%\noindent [17] L. De Pascale, M. S. Gelli, L. Granieri, Minimal measures, one-dimensional currents and the Monge-kantorovich problem, Calc. Var. and Partial Differential Equation, {\bf 27}, 363-388, 2006.

%\noindent [18] I. Ekeland, R. Temam, Convex analysis and variational problems, Amsterdam, 1976.

%\noindent [11] S. N. Ethier, T. G. Kurtz, Markov processes, Wiley, New York, 1986.

%\noindent [19] J. Feng, T. Nguyen, Hamilton-Jacobi equations in space of measures associated with a system of conservation laws, Journal de Math\'ematiques pures et Appliqu\'ees, {\bf 97}, 318-390, 2012. 

\noindent [14] I. Ekeland, Elements d'economie mathematique, Hermann, Paris, 1979.

\noindent [15] W. Gangbo, A. Swiech, Existence of a solution to an equation arising from the theory of mean field games, preprint 2014. 

\noindent [16] W. Gangbo, A. Swiech, Metric viscosity solutions of Hamilton-Jacobi equations, preprint 2014.

\noindent [17] W. Gangbo, A. Swiech, Optimal transport and large number of particles, Discrete and Continuous Dynamical Systems, {\bf 34,4}, 1387-1441, 2014.

\noindent [18] W. Gangbo, A. Tudorascu, Lagrangian dynamics on an infinite-dimensional torus; a weak KAM theorem, Adv. Math., {\bf 224}, 260-292, 2010.

\noindent [19] W. Gangbo, A. Tudorascu, Weak KAM theory on the Wasserstein torus with multi-dimensional underlying space, Comm. Pure Appl. Math., {\bf 67-3}, 408-463, 2014.

\noindent [20] Y. Giga, N. Hamamuki, A. Nakayasu, Eikonal equations in metric spaces, preprint.

\noindent [21] D. Gomes, L. Nurbekian, On the minimizers of variational problems in Hilbert spaces, Calc. Var., {\bf 52}, 65-93, 2014.

%\noindent [17] I. M. Gel'fand, A. M. Yaglom, Integration in functional spaces and its applications in Quantum Physics, J. Math. Phys. {\bf 1}, 48-69, 1960.

%\noindent [18] W. Gangbo, T. Nguyen, A. Tudorascu, Hamilton-Jacobi equations in the Wasserstein space, Methods Appl. Anal., {\bf 15}, 155-183, 2008.

%\noindent [19] D. Gomes, A stochastic analog of Aubry-Mather theory, Nonlinearity, {\bf 15}, 581-603, 2002.

%\noindent [21] D. Gomes, E. Valdinoci, Entropy penalization method for Hamilton-Jacobi equations, Advances in Mathematics, {\bf 215}, 94-152, 2007.

%\noindent [H] T. Hida, Brownian motion, Springer, Berlin, 1980.

%\noindent [22] R. Jordan, D. Kinderlehrer, F. Otto, The variational formulation of the Fokker-Planck equation, SIAM J. Math. Anal., {\bf 29}, 1-17, 1998.

%\noindent [22] E. Nelson, Dynamical theories of Brownian motion, Princeton, 1967.

\noindent [22] H. L. Royden, Real Analysis, China Machine Press, 2004. 

\noindent [23] C. Villani, Topics in optimal transpotation, Providence, R. I., 2003.

\end